\documentclass[preprint,12pt]{elsarticle}

\usepackage{amssymb}
\usepackage{amsmath}
\usepackage{subcaption}
\usepackage[dvipsnames]{xcolor}

\usepackage[normalem]{ulem}
\usepackage[utf8]{inputenc}

\journal{Wave Motion}

\begin{document}
\begin{frontmatter}


\title{Discretization in Multilayered Media with High Contrasts: Is It All About the Boundaries?}
\author[insa]{Camille Carvalho}
\author[poems]{Stéphanie Chaillat}
\author[ucm]{Elsie Cortes\corref{cor1}\fnref{fn1}}
\author[ucm]{Chrysoula Tsogka}

\cortext[cor1]{Corresponding author}
\fntext[fn1]{Email: ecortes7@ucmerced.edu}

\affiliation[insa]{
  organization={INSA Lyon, CNRS, Ecole Centrale de Lyon, Universite Claude Bernard Lyon 1, Université Jean Monnet, ICJ UMR5208},
  city={Villeurbanne},
  postcode = {69621},
  country = {France}
}
\affiliation[ucm]{
  organization={University of California, Merced},
  addressline={5200 Lake Rd},
  city={Merced},
  state={CA},
  postcode={95343},
  country={USA}
}

\affiliation[poems]{
  organization={Laboratoire POEMS, CNRS-ENSTA-INRIA, Institut Polytechnique de Paris},
  addressline={828 Bd des Marechaux},
  city={Palaiseau},
  postcode = {91120},
  country = {France}
}
\title{}

\begin{abstract}

Wave propagation in multilayered media with high material contrasts poses significant numerical challenges, as large variations in wavenumbers lead to strong reflections and complex transmission of the incoming wave field.  To address these difficulties, we employ a boundary integral formulation thereby avoiding volumetric discretization.  In this framework, the accuracy of the numerical solution depends strongly on how the material interfaces are discretized.  In this work, we demonstrate that standard meshing strategies based on resolving the maximum wavenumber across the domain become computationally inefficient in multilayered configurations, where high wavenumbers are confined to localized subdomains. Through a systematic study of multilayer transmission problems, we show that no simple discretization rule based on the maximum wavenumber or material contrasts emerges.  Instead, the wavenumber of the background (exterior) medium plays a dominant role in determining the optimal boundary resolution.  Building on these insights, we propose an adaptive approach that achieves uniform accuracy and efficient computation across multiple layers.  Numerical experiments for a range of multilayer configurations demonstrate the scalability and robustness of the proposed approach.

\end{abstract}


%
%
%
%
%


\end{frontmatter}



\section{Introduction}
\label{sec1}

We consider wave propagation in 2D multilayered media modeled by the Helmholtz equation.  The domain is modeled as a sequence of \(N + 1\) concentric layers, separated by smooth interfaces across which transmission conditions are imposed.  Such multilayer models are important in various applications such as noninvasive medical imaging \cite{bertolottiNoninvasiveImagingOpaque2012}\cite{eliseeInnovativeBoundaryIntegral2011}, non-destructive testing \cite{caiMethodologyEvaluatingDamage2024}\cite{caorsiMicrowaveImagingCylindrical2004}, and especially invisibility cloaking \cite{woodMetamaterialsInvisibility2009}\cite{pendryControllingElectromagneticFields2006}.  Our long-term motivation is related to optical cloaking, which can be obtained with metamaterials (artificially structured materials with deliberately engineered variations in material properties that enable complex wave behavior).   Simulating metamaterials and optimizing their design are of particular interest as an insight into their underlying physics, bridging the gap between theoretical understanding and practical implementation.  One of the challenges in designing metamaterials is dealing with large contrasts in wavenumbers $k_i$, \(i \in [0, N]\) across the interfaces (e.g., \(\frac{k_{i+1}}{k_i} \ll 1\) or \(\frac{k_{i+1}}{k_i} \gg 1\)), since they lead to strong reflection or complex transmission of the incident wave field.

Configurations with high material contrasts between layers inherently give rise to both high and low wavenumber regimes.  Mesh-based numerical approaches encounter significant challenges in the high-\(k_i\) regions, where the shorter wavelength \(\lambda_i\) of the wave field demands finer discretization and greater computational effort.  A common approach to discretizing wave problems is to maintain mesh size \(h\) such that \(kh \ll 1\). This conservative choice  prioritizes accuracy by ensuring sufficient resolution of the solution, even at the cost of overestimating the number of points per wavelength \cite{Ihlenburg1998}.  For multiple layers with different wavenumbers, the discretization is often governed by the maximum wavenumber \(k_{max}\) across all layers, i.e., according to the most constraining domain. This naive approach constrains the mesh to uniformly resolve the highest wavenumber although some parts of the domain correspond to lower wavenumbers. Furthermore, the pollution effect implies that the number of degrees of freedom must grow faster than the magnitude of \(k_i\) to maintain accuracy \cite{feng_sixth_2023}.  This effect implies that increasingly fine volumetric meshes are necessary to maintain accuracy.  Consequently, for multilayer problems, volumetric methods suffer from rapidly increasing computational costs since achieving accurate resolution requires mesh refinement even in regions where the wave field remains relatively simple.  While domain decomposition approaches combined with adaptive mesh refinement in each subdomain can provide effective solutions to this problem (see \cite{nataf2025} and references therein), we follow a different approach here. This deterioration in computational scalability motivates our use of a boundary integral formulation, where discretization is confined to material interfaces instead of the entire volume.

Boundary element methods (BEM) offer an efficient framework for approximating solutions to partial differential equations (PDE).
They rely on boundary integral equations (BIE), which reformulate the original PDE problem as an equation defined solely on the boundaries or interfaces.
BEM have been   applied to a wide range of applications e.g. \cite{Bonnet1999-bf}) including wave propagation problems, as they naturally enforce radiation conditions and largely avoid the pollution effect \cite{kress_boundary_1991,Galkowski2023}. A key advantage of this approach is the reduction of dimensionality: only boundary data must be discretized, avoiding the need to mesh the entire computational domain, unlike   Finite Element Methods (FEM) or Finite-Difference (FD) methods.  This makes BEM especially attractive for layered media  \cite{chen_thermal_2001,liu_bem_2019}, although multilayer configurations   lead to   densely populated matrices \cite{magoules_boundary_2008}.  Moreover, the singular or nearly-singular kernels appearing in the BIE require careful treatment \cite{barnett_evaluation_2014}. 
While techniques for weakly \cite{kress_boundary_1991} and nearly-singular kernels \cite{Klckner2013,Carvalho2020} are well-established for single-domain problems, multilayered domains require a boundary integral framework configured to accurately handle multiple interfaces in relation to each other. We have previously derived a modified boundary equation method to find accurate solutions to the multilayered transmission problem when interfaces are close to each other \cite{cortesFastAccurateBoundary2024}.  By reducing the influence of nearly singular kernels on the boundary solution, this method avoids the need for fine interface refinement to obtain a precise solution.

One of the critical factors in BEM is the discretization of the interfaces, which strongly influences both the accuracy of the solution and the efficiency of the computation. Coarse meshes or meshes that fail to capture the boundary geometry or key solution features can produce significant errors, while overly refined meshes increase computational cost. In the current work, we tackle this challenge with adaptive mesh refinement.  We propose an effective, systematic way to optimize the discretization by refining the mesh only where it is needed, improving both the quality of our solutions and the computational performance.  In general, mesh adaptation in the context of BEM remains an open area of research where several theoretical and practical questions are being explored \cite{wartmanFastEEGMEG2025,wartmanAdaptiveHrefinementMethod2024}.  While no universal method exists in the BEM community, the situation is even more challenging for multilayered domains. In this work, we specifically investigate the use of metric-based mesh refinement for this problem \cite{mirebeau2011,chaillat_metric-based_2018}. 

By focusing on the discretization of the \textbf{boundaries} rather than the entire domain, this work presents a new perspective on discretization of multilayer media for wave propagation problems.  In Section \ref{sec2},  we begin by formulating the multilayer transmission problem with boundary integral equations   and present high-order quadrature rules for their discretization.  In Section \ref{sec:bie-based mesh optimization}, we explore boundary mesh optimization, using reference analytic solutions for single-material domains to construct \emph{a posteriori }estimates that provide insight into solution accuracy for multilayer configurations with increasing contrasts.  We observe that the accuracy of the numerical solutions is not driven by the highest wavenumber. Building on this, in Section \ref{sec:adaptive mesh generation}, we extend our focus to adaptive mesh generation,  automatically refining boundary meshes in configurations with multiple boundaries and more complex geometries.  

Section 5 concludes with a discussion of possible directions for future research.

\section{The Multilayer Transmission Problem}
\label{sec2}
\subsection{Problem Formulation for Concentric Regions}

\begin{figure}[!h]
      \centering
      \includegraphics[width=3in]{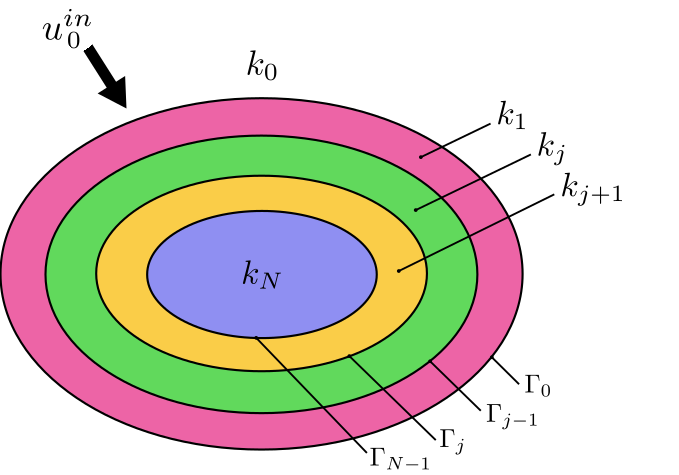}
      \caption{Schematic of the concentric layered transmission problem.  Each interface $\Gamma_j$ separates regions with wavenumbers $k_j$ and $k_{j+1}$. The exterior field $u_0$ includes an incoming incident wave $u_0^{in}$.}
      \label{fig:transmission problem}
\end{figure}

Consider a domain composed of \(N + 1\) concentric regions in \(\mathbb{R}^2\), separated by closed, smooth interfaces \(\Gamma_j\) for \(j = 0, \hdots, N - 1\) as illustrated in Figure \ref{fig:transmission problem}.  Let \(\Omega_j \subset\mathbb{R}^2\) for \(j = 0, \hdots, N\) denote the concentric subdomains, such that \(\Omega_0\) is the outermost region, or the background medium, and \(\Omega_N\) is the innermost region.  The interfaces \(\Gamma_j\) separate \(\Omega_j\) and \(\Omega_{j+1}\) for \(j = 0, \hdots, N - 1\).  Each region is characterized by a wavenumber \(k_j \in \mathbb{R}_+\) (determined by its material properties).  We define the \textbf{contrast} between the material properties of adjacent regions \(\Omega_{j}\) and \(\Omega_{j+1}\) as \(\beta_{j} := k^2_{j+1} / k^2_{j}\).

Working in the time-harmonic regime, the total field solution in \(\Omega_0\) is the sum of the incident and scattered fields, and the incident field is a plane wave of the form \(u_0^{\mathrm{in}} = e^{\mathrm{i} k_0 \mathbf{a} \cdot \mathbf{x}}\) (\footnote{In this work, we choose \(\mathbf{a} = (0,1)\) so that the incident field \(u^{\mathrm{in}} = e^{\mathrm{i} k_0 y}\) for \(\mathbf{x} := (x, y)\) is an incoming plane wave propagating in the positive y-direction. But all results in this paper hold for any \(\mathbf{a}\).}).
In each layer \(\Omega_j\) the total field \(u_j\) satisfies the Helmholtz equation \cite{colton_helmholtz_2019}:

\begin{equation}
     \Delta u_j + k^2_j u_j = 0 \quad \text{in} \quad \Omega_j.
     \label{eq:Helmholtz Eq}
\end{equation}
Furthermore, across each interface we impose the following transmission conditions:
\begin{equation}
\begin{aligned}
u_{j+1} = u_{j }  \quad & \text{on} \quad \Gamma_{j} & \quad j = 0, ..., N-1, \\
  \frac{\partial u_{j+1}}{\partial n_{j}} = \beta_{j } \frac{\partial u_{j }}{\partial n_{j }}  \quad &  \text{on} \quad \Gamma_{j} &\quad  j = 0, ..., N-1,
\end{aligned}
\label{eq:transmission problem}
\end{equation}
where \(n_j\) denotes the outward unit normal on \(\Gamma_j\).  Additionally, we assume the scattered solution in $\Omega_0$, $u_0^{\mathrm{sc}} = u_0 - u_0^{\mathrm{in}}$, satisfies the Sommerfeld radiation condition:

\begin{equation}\lim_{r \to \infty} \sqrt{r} \left( \frac{\partial u_0^{\mathrm{sc}}}{\partial r} - ik_0u_0^{\mathrm{sc}} \right) = 0  .
\label{eq:sommerfeld}
\end{equation}
This condition ensures that the outgoing wave solution decays as \(r\) tends to infinity.

With the multilayer transmission model in place, we now turn to finding the boundary integral formulation for the BEM.

\subsection{Boundary Integral Equation Solutions}

BIE, as the name implies, are representations of the solutions in terms of boundary integrals.  For the total field satisfying the Helmholtz equation in each layer, a boundary integral formulation may be obtained using a combination of single- and double-layer boundary integral operators \(\mathcal{S}_{i,j}\) and \(\mathcal{D}_{i,j}\) respectively \cite{kress_boundary_1991}.  They are defined as 

\begin{align}
    \mathcal{S}_{i,j}[\mu](\mathbf{x}) &= \int_{\Gamma_i} \Phi_j(\mathbf{x},\mathbf{y}) \mu(\mathbf{y}) d\Gamma_{i,\mathbf{y}} , &\quad \mathbf{x} \in \mathring{\Omega}_j,\\
        \mathcal{D}_{i,j}[\mu](\mathbf{x}) &= \int_{\Gamma_i} \frac{\partial \Phi_j}{\partial n_i(\mathbf{y})}(\mathbf{x}, \mathbf{y}) \mu(\mathbf{y}) d\Gamma_{i,\mathbf{y}} ,&\quad \mathbf{x} \in \mathring{\Omega}_j ,
    \label{eq:layer potentials}
\end{align}

where  \(\Gamma_i\) (for \(i = j - 1\)  or \(j\)) denotes the interface on which they are defined, and \(\Phi_j\) is the fundamental solution to the Helmholtz equation in \(\Omega_j\). In 2D the fundamental solution is given by:
\begin{equation}
    \Phi_j(\mathbf{x},\mathbf{y}) = \frac{\mathrm{i}}{4}H_0^{(1)}(k_j|\mathbf{x}-\mathbf{y}|)
    \label{eq:fundamental solution}
\end{equation}
where $H_0^{(1)}$ denotes the Hankel function of first kind.
There are multiple boundary integral formulations, but for the one we choose the unknown densities  $\mu \in  C^2(\Gamma_i)$ in \(\mathcal{S}_{i,j}\) and \(\mathcal{D}_{i,j}\)  define the boundary data.   More precisely, in the double-layer integral operator, $\mu$ corresponds to the trace $u_j$, and in the single-layer integral operator it corresponds to the normal trace \(\tfrac{\partial u_j}{\partial n_i}\).
For simplicity, we denote the traces \(u_j(\mathbf{y})\) and \(\tfrac{\partial u_j}{\partial n_i}(\mathbf{y})\) respectively by
\begin{equation*}
u_{i, j} := u_j(\mathbf{x}), \quad \text{and}\quad  \partial_n u_{i, j} := \frac{\partial u_j}{\partial n_i}(\mathbf{x}), \quad \quad \mathbf{x} \in \Gamma_i.
\end{equation*}

Using the representation formula \cite{kress_boundary_1991}, and utilizing the layer potentials, we construct the following boundary integral representations for the solutions for \(j = 1, \hdots, N - 1\):
\begin{equation}
\resizebox{0.9\textwidth}{!}{$
\begin{aligned}
    u_0(\mathbf{x}) &= u_0^{in}(\mathbf{x}) + \mathcal{D}_{0,0}[u_{0,0}](\mathbf{x}) + \mathcal{S}_{0,0}[\partial_n u_{0,0}](\mathbf{x}) & \text{in } \mathring{\Omega}_0, \\
    u_j(\mathbf{x}) &= -\mathcal{D}_{j-1,j}[u_{j-1,j}](\mathbf{x}) + \mathcal{S}_{j-1,j}[\partial_n u_{j-1,j}](\mathbf{x}) \\ 
    &\quad + \mathcal{D}_{j,j}[u_{j,j}](\mathbf{x}) - \mathcal{S}_{j,j}[\partial_n u_{j,j}](\mathbf{x}) & \text{in } \mathring{\Omega}_j, \\
    u_N(\mathbf{x}) &= -\mathcal{D}_{N-1,N}[u_{N-1,N}](\mathbf{x}) + \mathcal{S}_{N-1,N}[\partial_n u_{N-1,N}](\mathbf{x}) & \text{in } \mathring{\Omega}_N
\end{aligned}
$}
\label{eq:bie_layers}
\end{equation}

The solution in each domain is determined from the traces of \(u\) and \(\partial u\). To obtain these traces \(u_{i,j}\) and \(\partial u_{i,j}\)  we construct a system of representation equations defined on the interfaces \(\Gamma_j\) in terms of layer potentials.  The layer potentials can be obtained by taking the limits of the single- and double-layer layer integral operators on each side of the boundaries. The single-layer integral operator is a continuous operator such that  

\begin{equation*}\mathcal{S}_{i,j}[\mu](\mathbf{x}) \underset{\mathbf{x}\rightarrow \mathbf{x}^{b,\pm}}{\longrightarrow} S_{i,j}[\mu](\mathbf{x}^b),\end{equation*}

but the double-layer operator exhibits a jump such that

\begin{equation*}\mathcal{D}_{i,j}[\mu](\mathbf{x}) \underset{\mathbf{x}\rightarrow \mathbf{x}^{b, \pm}}{\longrightarrow} \pm \frac{1}{2}\mu(\mathbf{x}^b) + D_{i,j}[\mu](\mathbf{x}^b)\end{equation*}

where $\mathbf{x}^b$ is on the boundary. 

By introducing the transmission conditions \eqref{eq:transmission problem} in the boundary integral equation system \eqref{eq:boundaryreps}, we obtain a solvable system of \(2N\) equations with \(2N\) unknowns, namely the traces \(u_{i,j}\) and normal traces \(\partial u_{i,j}\) respectively that we can directly solve for these quantities.  For example, if we substitute \(u_{j-1,j} = u_{j-1,j-1}\) and \(\partial_n u_{j-1,j} = \beta_{j-1} \partial_n u_{j-1,j-1}\), we obtain the following boundary integral equations:
\begin{equation}
\resizebox{0.9\textwidth}{!}{$
\begin{aligned}
    \tfrac{1}{2}u_0(\mathbf{x}^b) &= u_0^{in}(\mathbf{x}^b) + D_{0,0}[u_{0,0}](\mathbf{x}^b) + S_{0,0}[\partial_n u_{0,0}](\mathbf{x}^b) & \text{on } \Gamma_0, \\
    \tfrac{1}{2}u_j(\mathbf{x}^b) &= -D_{j-1,j}[u_{j-1,j-1}](\mathbf{x}^b) + \beta_{j-1}S_{j-1,j}[\partial_n u_{j-1,j-1}](\mathbf{x}^b) \\
    &\quad + \mathcal{D}_{j,j}[u_{j,j}](\mathbf{x}^b) - \mathcal{S}_{j,j}[\partial_n u_{j,j}](\mathbf{x}^b) & \text{on } \Gamma_{j-1}, \\
    \tfrac{1}{2}u_j(\mathbf{x}^b) &= -\mathcal{D}_{j-1,j}[u_{j-1,j-1}](\mathbf{x}^b) + \beta_{j-1}\mathcal{S}_{j-1,j}[\partial_n u_{j-1,j-1}](\mathbf{x}^b) \\
    &\quad + D_{j,j}[u_{j,j}](\mathbf{x}^b) - S_{j,j}[\partial_n u_{j,j}](\mathbf{x}^b) & \text{on } \Gamma_j, \\
    \tfrac{1}{2}u_N(\mathbf{x}^b) &= -D_{N-1,N}[u_{N-1,N}](\mathbf{x}^b) + \beta_{N-1} S_{N-1,N}[\partial_n u_{N-1,N-1}](\mathbf{x}^b) & \text{on } \Gamma_{N-1}
\end{aligned}
$}
\label{eq:boundaryreps}
\end{equation}

There exist various approaches to discretizing this system.  Since we consider multilayer problems defined by smooth periodic boundaries, we choose to parameterize the boundary integrals in terms of \(\mathbf x(\theta) = (x_1(\theta), x_2(\theta))\) for \(0 \leq \theta \leq 2\pi\).  A well-known approach for the defined boundary discretization is proposed by Kress \cite{kress_boundary_1991}.  It consists of expressing the boundary integrals in terms of trigonometric polynomials equi-spaced on \(\theta_m\) such that
\begin{equation}
    \theta_m := \frac{2\pi m }{M} \quad \text{ for } m = 0, \hdots, M - 1.
    \label{eq:theta distribution}
\end{equation}

A key point for BEM is the evaluation of the (singular) boundary integrals.  In the case of 2D Helmholtz, layer potentials in \eqref{eq:boundaryreps} contain a kernel with a logarithmic singularity.  An effective treatment method is the Kress Quadrature Rule (KQR) \cite{kress_boundary_1991}.  KQR splits each singular kernel into an analytically integrable part plus a smooth remainder, allowing the use of corrected trapezoidal weights on the latter. 

From an algebraic point of view, the size of the system of equations to solve is determined by the number of boundaries \(N\) and the number of equispaced points \(M\) on each interface.  Once we determine the unknown traces and normal traces, the solution is obtained using the boundary integral representation in the volume \eqref{eq:bie_layers}, approximated with the periodic trapezoid rule (PTR) where we evaluate on the same uniform distribution of \(\theta_m\) in \eqref{eq:theta distribution}:

\begin{equation}
\int_0^{2\pi} f(\theta) d\theta \approx \Delta \theta \sum_{m = 0}^{M} f(\theta_m) = \frac{2\pi}{M} \sum_{m = 0}^{M - 1} f(\theta_m)
\label{eq:ptr approx}
\end{equation}

For smooth, closed boundaries, the approach gives spectral convergence with respect to the number of quadrature points \(M.\)

A key question for this problem is how to determine the optimal number of points \(M\) to use on each interface.  Traditionally, a general rule based on the maximum wavenumber \(k_{max}\) in the domain is used to set M, either uniformly across all interfaces or scaled according to the size of each boundary.  When \(k_{max}\) is large, the field becomes highly oscillatory, requiring a greater number of points to be accurately represented.  One of our objectives is to go beyond the standard approach and determine the optimal number of points \(M_i\) on each interface \(\Gamma_i\), where optimal is defined as the minimal \(M_i\) that achieves a desired order of accuracy in the boundary data at that interface.  In particular, we aim to investigate whether a rule that accounts for the contrast at each interface, \(\beta_j :=\frac{k_{j+1}^2}{k_j^2}\), can be inferred.  This is a largely unexplored area for wave propagation problems in multilayer media.   

In the next section, we introduce an approach to identify the optimal \(M_i\) for various contrasts, including comparisons to the sound-hard case, in order to examine any trends observed in multilayer problems.

\section{Boundary Mesh Optimization for BEM}
\label{sec:bie-based mesh optimization}
\subsection{A Preliminary Approach Based on the Reference Solution}

We focus on configurations with circular boundaries, for which analytic solutions can be derived (see Section~\ref{sec:analytic solutions}) and used as reference solutions to gain further insight into the problem. We aim to determine the minimum number of points \(M_i\) on each interface \(\Gamma_i\) required to achieve a target accuracy. Here, accuracy is measured as the absolute \(L^2\) error between the BIE solution of \eqref{eq:boundaryreps} on the boundaries and the corresponding analytic solution \(u^{ana}\) \eqref{eq:analytic solutions}.  Our naive preliminary approach to determine the optimal discretization consists of iteratively increasing \(M_i\) until the inequality
\begin{equation}
    ||u_j - u_j^{ana}||_{L^2({\Gamma_i)}} \leq \epsilon
    \label{eq:error_ineq}
\end{equation}
is satisfied, with \(\epsilon\) the desired accuracy (between \(10^{-6}\) and \(10^{-8}\) in this work). We stop the iterative process when \eqref{eq:error_ineq} is satisfied for all boundaries.

In the context of multilayer transmission problems, evaluating errors on the interfaces is essential for assessing how accurately the discretization captures the complex behavior of the solution in the surrounding volume.
To evaluate the absolute \(L^2\) error  \eqref{eq:error_ineq}, we use the explicit parameterization of the interfaces \(\Gamma_i\) for \(i = j-1\) or \(j\).  It follows that
\begin{equation}
    ||u_j - u_j^{ana}||_{L^2(\Gamma_i)} = \left( \int_{\Gamma_i} |u_j(\mathbf{x}_j(\theta)) - u_j^{ana}(R_i, \theta) |^2 d\theta \right)^{1/2}
    \label{eq:parameterized error}
\end{equation}
for \(\mathbf{x}_j(\theta):= (R_i \cos\theta, R_i \sin\theta) \).

Since the interfaces $\Gamma_i$ are periodic and parameterized in terms of $\theta$ over the interval $[0, 2\pi)$, we use PTR to evaluate the integral in \eqref{eq:parameterized error}.  The boundary is discretized using $M_i$ equispaced points \(\mathbf x_j(\theta_m)\) for \(m = 0, \hdots, M_i\) and \(\theta_m\) defined as in \eqref{eq:theta distribution}, so that \eqref{eq:parameterized error}
\begin{equation}
||u_j - u_j^{ana}||_{L^2(\Gamma_i)} \approx \left( \frac{2\pi}{M_i} \sum_{m = 0}^{M_i - 1} |u_j(\mathbf{x}_j(\theta_m)) - u_j^{ana}(R_i, \theta_m) |^2\right)^{1/2}.
\label{eq:error_approx}
\end{equation}

\subsection{Observations from Analytic Solutions}

\begin{figure}[!b]
    \centering   
    \hfill
    \begin{subfigure}[t]{0.49\textwidth}
        \centering
        \includegraphics[width=0.65\textwidth]{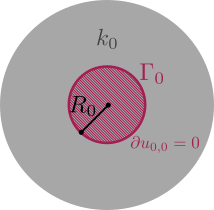}
        \caption*{2(a): Sound-hard boundary configuration}
    \end{subfigure}
    \begin{subfigure}[t]{0.5\textwidth}
        \centering
        \includegraphics[width=0.67\textwidth]{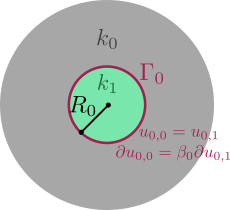}
        \caption*{2(b): Transmission boundary configuration}
    \end{subfigure}
    \caption{Schematics of the scattering domain with radius $R_0 = 4$, showing two types of boundary conditions:  (a) a sound-hard obstacle and (b) a penetrable interface with transmission conditions.}
    \label{fig:problem schema}
\end{figure}

We consider both the sound-hard boundary case imposing \(\partial_n u_{0,0}= 0\) on \(\Gamma_0\) and the transmission boundary case (see \ref{sec:analytic solutions} and \ref{app:analytic solutions.} for expressions of the associated analytic solutions).  The goal of this comparison is to reveal how the properties of an additional layer affect the solution in the previous layer.  Rather than focusing solely on the layer with the highest wavenumber, where oscillations are strongest, this analysis helps to understand how the solution evolves across multiple layers.

We consider a simple case with only one interface (a circular boundary at \(R_0 = 4\)) defining two materials. We analyze the accuracy for different contrasts,  $\beta_0 = k_1^2 / k_0^2$ , where $k_0$ and $k_1$ are the wavenumbers in the outer and inner regions respectively as shown in Figure \ref{fig:problem schema}(b).  We examine all the cases in Table \ref{tab:1bdycases}.  In particular, we explore cases obtained by swapping the interior and exterior wavenumbers.  The aim is to examine how the contrast influences the wave propagation and scattering behaviors. All cases are compared to the corresponding sound-hard case, using their respective $k_0$ in the exterior region (see Figure \ref{fig:problem schema}(a)).

\begin{table}[!t]
\centering
\begin{tabular}{|c c c c|}
\hline
Case & $\mathbf{k} = (k_0, k_1,...,k_N)$ & $\mathbf{r} = (R_0, ..., R_{N-1})$& $\mathbf{b} = (\beta_0, ..., \beta_{N-1})$\\
\hline
1 & $(2,6)$  & $4$ & 9\\ 
2 & $(6,2)$  & $4$ & $\tfrac{1}{9}$\\
3& $(2,10)$ & $4$ & $25$ \\
4& $(10,2)$ & $4$ & $\tfrac{1}{25}$ \\
5& $(8, 2, 6)$ & $(6, 2)$ & $(\tfrac{1}{16},  9)$\\
\hline
\end{tabular}
\caption{Summary of single-boundary cases used in numerical experiments, with wavenumbers $\mathbf{k}$, boundary radius $R_0$, and contrast $\beta_0$.}
\label{tab:1bdycases}
\end{table}

\begin{figure}[!h]
    \centering
    \includegraphics[width=\textwidth]{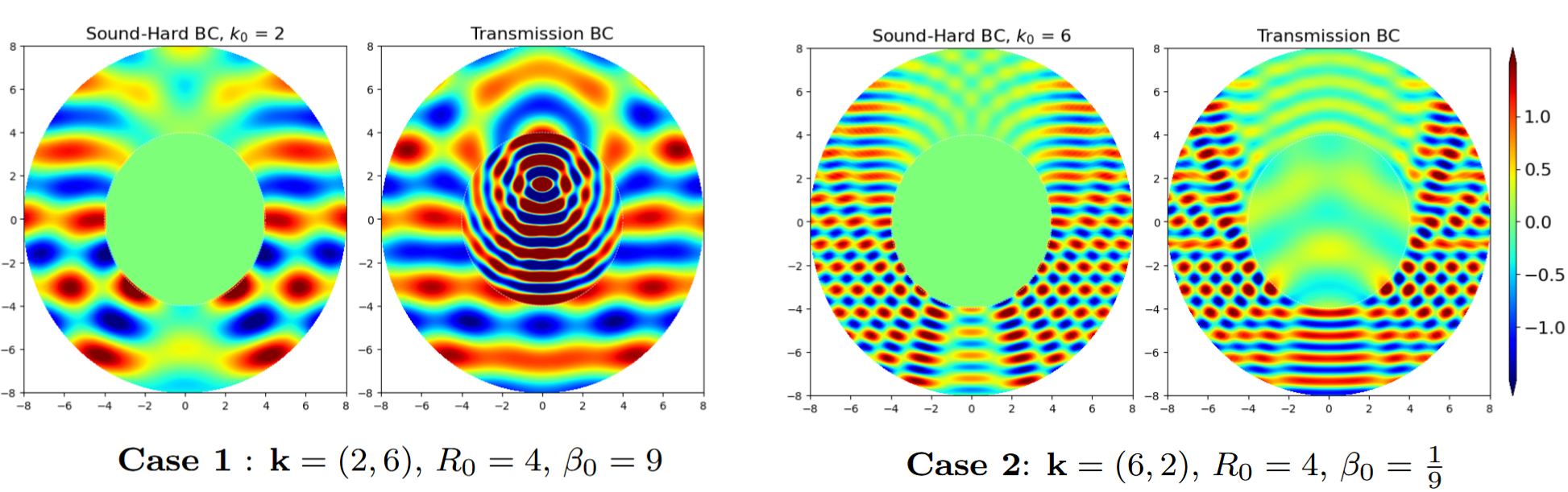}
    
    \caption{Real part of the total field analytic solutions for two contrast cases with $R_0 = 4$. The left plot on each sub-figure shows the solution for the sound-hard boundary condition (Neumann), while the right plot of each sub-figure shows the solution for the transmission boundary condition.}
    \label{fig:analytic sets 1}
\end{figure}
In Figure \ref{fig:analytic sets 1}, we plot the corresponding total field solutions for \textbf{Cases 1 and 2}.  In \textbf{Case 1}, we observe that for the transmission case the solution exhibits strong oscillations within the interior region, although much of this complexity is constrained  by the boundary and the lower wavenumber of the exterior region. Although this suggests that finer discretization may be necessary, the oscillatory behavior along the boundary remains comparable to the reference solution for the sound-hard boundary condition.  In contrast, in \textbf{Case 2}, the transmission boundary condition produces wave fields with oscillations comparable to those of the exterior for the sound-hard case.  This supports the expectation that resolution should be governed by the higher wavenumber. However, the decrease in oscillations near the boundary indicates a more subtle effect.

\begin{figure}[!b]
    \includegraphics[width=\textwidth]{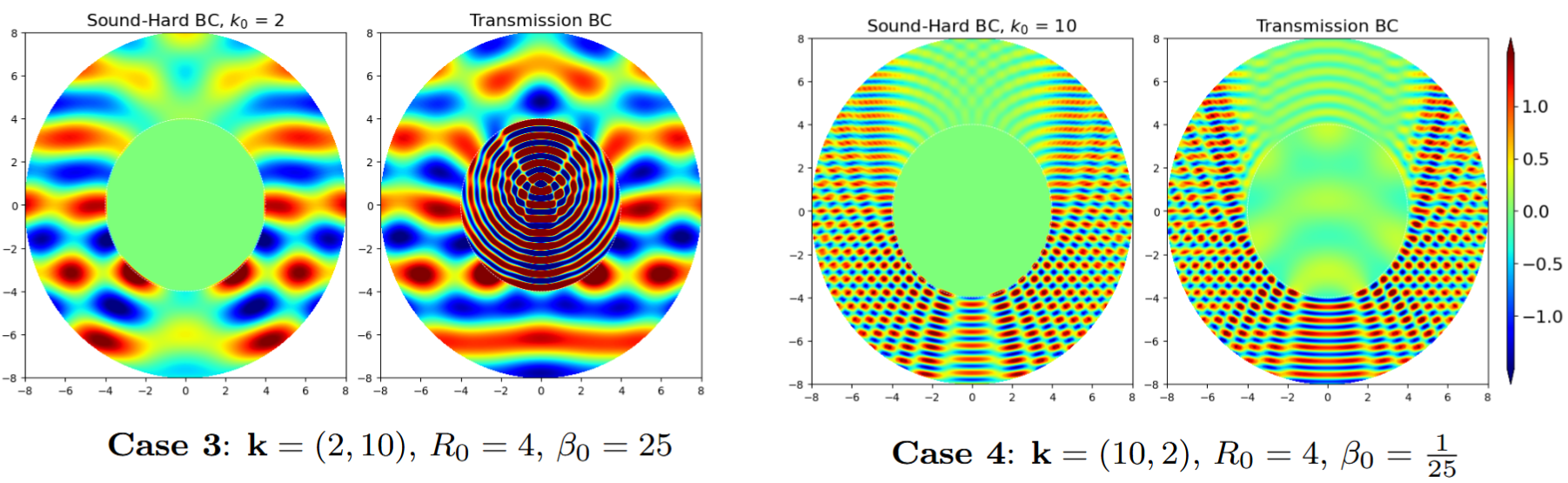}
    \caption{Real part of the total field analytic solutions for a higher contrast case with sharper contrasts $\beta_0$ as given.}
    \label{fig:analytic sets 2}
\end{figure}

For the sharper contrasts in \textbf{Cases 3 and 4}, shown in Figure \ref{fig:analytic sets 2}, we draw similar conclusions.  However,  for subplots with the transmission boundary condition, much of the solution behavior is preserved in regions where \(k=2\) despite the sharp interface contrast.  For example, comparing the solution in \(\Omega_0\) for \textbf{Case 1} and \textbf{Case 3} reveals little to no visual difference.  This suggests that the level of discretization required to achieve the same accuracy on \(\Gamma_0\) for \(\beta_0=9\) and \(\beta_0=25\) is more similar than one might expect given the difference in contrast.  To determine whether this is indeed the case, we examine the optimal number of boundary points required to achieve a fixed target accuracy.

\subsection{Observations from Boundary Optimization}

Now we check the number of points on the boundary required to achieve an accuracy of \(10^{-6}\) using solutions approximated with BEM. Beginning with \textbf{Case 1} with the sound-hard boundary condition, we determine the required number of discretization points on the boundary, denoted \(M^{SH}_0\).  \(M_0^{SH}\) is obtained by incrementally refining the boundary discretization until the target \(L^2\) error, in \eqref{eq:error_approx}, falls below \(10^{-6}\).  Our tests indicate that a minimum of \(M_0^{SH} = 46\) points is required to reach this accuracy.  To estimate the number of points per wavelength \(N_{ppw}\) a discretization on the boundary \(\Gamma_j\) corresponds to, we use the relation

\begin{equation}
    M_j = N_{ppw} \cdot R_j \cdot \max(k_j) \quad \text{on } \Gamma_j.
    \label{eq:general rule}
\end{equation}

We need to determine $N_{ppw}$ for the target accuracy.  We obtain an estimate of it for the sound-hard case using:

\begin{equation*}
    M_0^{SH} \approx N_{ppw} \cdot R_0 \cdot k_0,
\end{equation*}
which gives \(N_{ppw} \approx 6\).  While this value is derived for a specific case, we use it as a practical reference point.  The goal is not to impose a precise standard for all cases, but to guide expectations for the discretization of the boundary.  Applying this to \textbf{Case 1} for transmission, when the highest wavenumber is \(k_1 = 6\), we expect to need approximately \(M_0 = 144\) points to achieve a similar accuracy.  However, the optimal value, as determined by iterating through increasing values of \(M\), is \( \boldsymbol{M^{tar}_0 = 84}\).

For \textbf{Case 2}, the optimal interface discretization for the transmission boundary condition is \(\boldsymbol{M^{tar}_0 = 108}\).  This value is higher than that obtained in the previous case, which involved a low-to-high wavenumber contrast.  This indicates that the exterior wavenumber influences the solution on both sides of the interface, so a higher exterior wavenumber leads to higher resolution demand (and the opposite is also true).   However, this result is also lower than what would be predicted by the common guideline for wave propagation problems based only on the highest wavenumber.  From \(M_0^{SH}=112\) we estimate \(N_{ppw} \approx 5\), which yields that the approximate number of points we should expect to obtain a similar accuracy for the transmission case is \(M_0 = 120\).  Adopting a more conservative baseline for the same pair of wavenumbers \((2,6)\) would suggest using our previously obtained \(M_0 = 144\).   Regardless of which assumption we follow, the observed optimal resolution falls below the estimates.  The general rule systematically overestimates the required boundary resolution in the multilayer setting.

\begin{figure}[!t]
      \centering
      \includegraphics[width=2.2in]{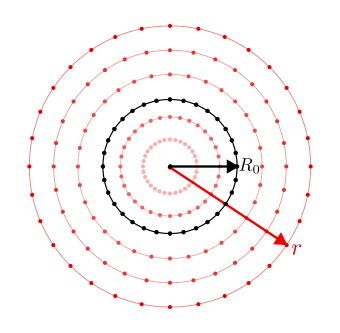}
      \caption{Schematic of the samplings for \(||u_j - u^{ana}_j||_{L^2(r)}\).  The boundary solution is computed on the black ring at radius \(R_0\), and the error between the off-boundary PTR and analytic solutions is evaluated at each sampling radii (red rings) separately.}
      \label{fig:radialerrorschema}
\end{figure}

To illustrate how accuracy on the boundary impacts the accuracy in the volume, we also determine the solution in the layers and compare the resulting error for different evaluation distances.  To do so, we determine the \(L^2\) error of the numerical solutions in the layers at various radial distances at radii r, both inward and outward from the boundary \(\Gamma_i\) (see Figure \ref{fig:radialerrorschema}).  The boundary solution is evaluated at the radii \(R_i\) using the BIE.  Then the fields in each layer $u_j$ are evaluated on the boundary-fitted grid using \eqref{eq:ptr approx}.  At each sampling radius for fixed values of \(r\), the numerical solution is compared to the analytic solution \eqref{eq:analytic solutions} on \(M_i\) values of \(\theta_m\). Specifically, the radial \(L^2\) error is computed using \eqref{eq:error_approx}, 
\begin{equation*}
||u_j - u_j^{ana}||_{L^2(r)} \approx \left( \frac{2\pi}{M_i} \sum_{m = 0}^{M_i - 1} |u_j(\mathbf{x}_j(\theta_m)) - u_j^{ana}(r, \theta_m) |^2\right)^{1/2}
\end{equation*}
for $i = j$, or otherwise $i = N - 1$ if within the region $\Omega_N$, which is not defined by an interior boundary.  Here, $\mathbf{x}_j(\theta_m) := (r\cos \theta_m, r\sin\theta_m)$.

\begin{figure}[!t]
    \centering
        \includegraphics[width=0.7\textwidth]{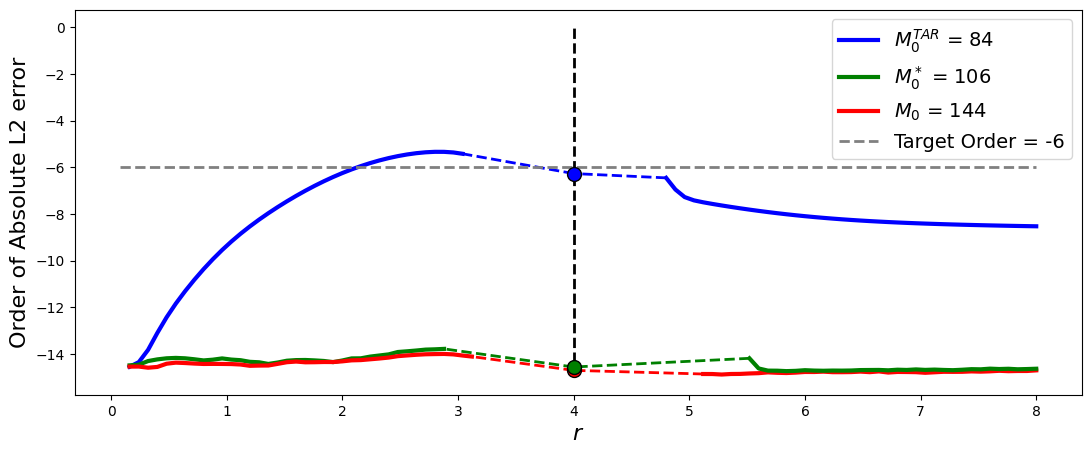}
        \caption*{Radial errors for \textbf{Case 1} : \(\mathbf{k} = (2,6)\), \(R_0  = 4\), \(\beta_0 = 9\)}
  
    \caption{Semilog plot of the absolute $L_2$ error of solutions over radii r.  The vertical dashed black line marks $R_0$, while the horizontal dashed gray line indicates the target order of error. Dashed segments denote skipped values due to close evaluation \cite{barnett_evaluation_2014}.}
    \label{fig:radial errors set 1}
\end{figure}

In Figure \ref{fig:radial errors set 1}, we examine the radial error behavior for \textbf{Case 1}. This figure shows that our BEM framework delivers accuracy in the volume that is comparable to, or even better than, that obtained on the boundary. This behavior comes from the smoothing effect of the boundary integral representation, combined with the spectral accuracy of the PTR for smooth and periodic interfaces, resulting in highly accurate solutions across the domain \cite{kress_boundary_1991}.  Using the estimate \(M_0 = 144\), we achieve machine precision at \(10^{-14}\) throughout the domain, rather than the intended target accuracy achieved with the optimal \(M^{tar}_0 = 84\).  While obtaining machine precision is not undesirable, this result indicates that the estimate over-resolves the solution.  The optimal number of points to reach machine precision is \(M^*_0 = 106\), which lies closer to the discretization needed for a lower order accuracy than our estimate \(M_0\). This example clearly demonstrates that only accounting for \(k_{max}\) can lead to unnecessary levels of discretization for multilayered media.

\begin{figure}[!h]
    \centering
        \includegraphics[width=0.7\textwidth]{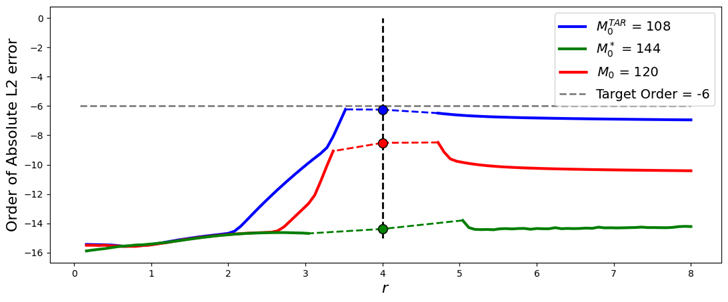}
        \caption*{Radial errors for \textbf{Case 2} : \(\mathbf{k} = (6,2)\), \(R_0  = 4\), \(\beta_0 = \frac{1}{9}\)}
  
    \caption{Semilog plot of the absolute $L_2$ error of solutions over radii r.}
    \label{fig:radial errors set 2}
\end{figure}

The radial errors for \textbf{Case 2} are shown in Figure \ref{fig:radial errors set 2}.  Estimating \(M_0\) from the sound-hard boundary discretization using \(N_{ppw} \approx 5\) results in a solution that falls below the target error threshold. A more conservative choice of \(N_{ppw} = 6\) (which is the average \(N_{ppw}\) estimated for every case in Table \ref{tab:1bdycases}) yields the optimal number of boundary points required for machine precision, \(M^*_0=144\).  In either case, the discretization remains overestimated, even though the solution behavior observed in Figure \ref{fig:analytic sets 1} shows little visual difference between sound-hard and transmission boundary conditions for Case 2.

\begin{figure}[!h]
    \centering
    \begin{subfigure}{\textwidth}
        \centering
        \includegraphics[width=0.7\textwidth]{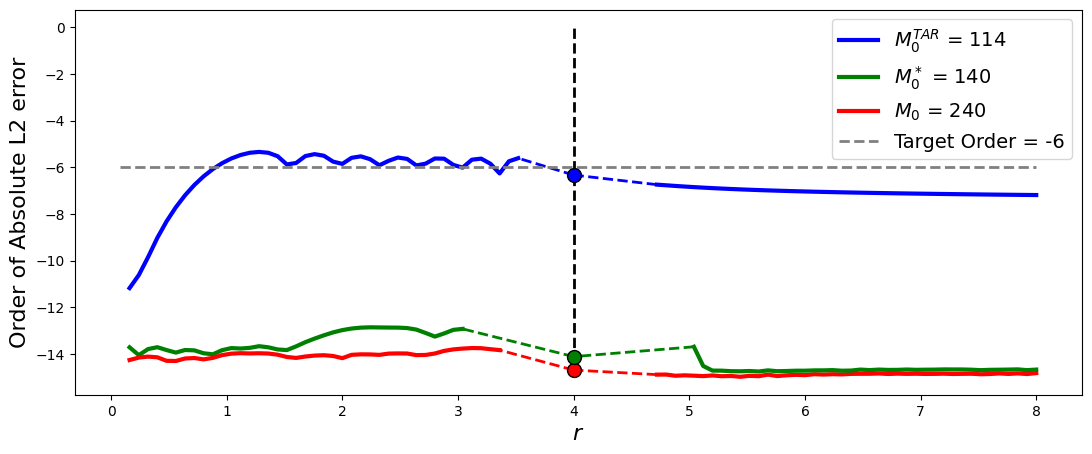}
        \caption*{Radial errors for \textbf{Case 3} : \(\mathbf{k} = (2,10)\), \(R_0  = 4\), \(\beta_0 = 25\)}
    \end{subfigure} 
   \begin{subfigure}{\textwidth}
        \centering
        \includegraphics[width=0.7\textwidth]{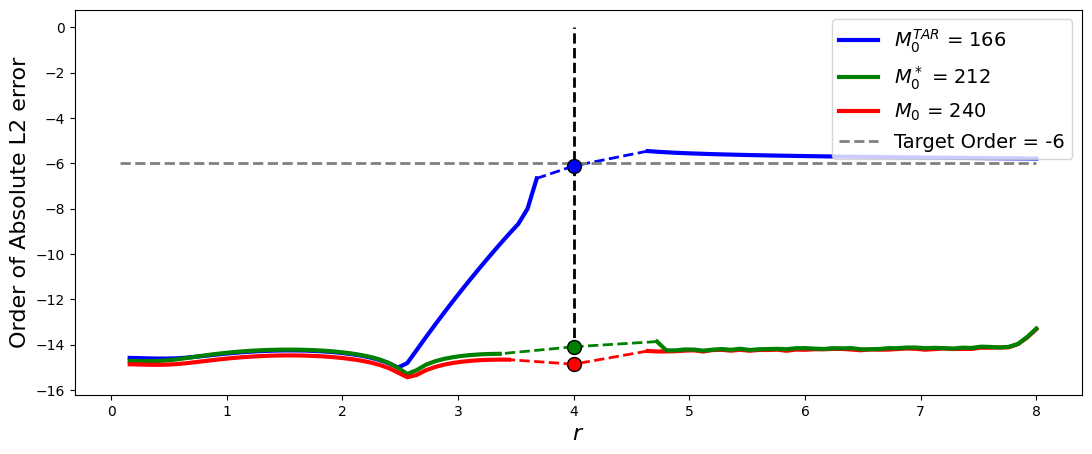}
        \caption*{Radial errors for \textbf{Case 4} : \(\mathbf{k} = (10, 2)\), \(R_0  = 4\), \(\beta_0 = 25\)}
    \end{subfigure}
    
    \caption{Semilog plot of the absolute $L_2$ error of solutions for the given cases computed at different radii r.}
    \label{fig:radial errors set 3}
\end{figure}

To see whether these observations are more pronounced for sharper contrasts, we also plot the radial errors for \textbf{Case 3} and \textbf{Case 4}.  In Figure \ref{fig:radial errors set 3}, \(M^{tar}_0\) in both cases deviates even further from the value obtained by the general rule, \(M_0\), in \eqref{eq:general rule}.  Even when targeting machine precision, the discretization is substantially overestimated relative to the true optimal value \(M^*_0\).  For example, in \textbf{Case 3}, \(M_0 = 240\) is nearly twice \(M^*_0 = 140\). These results indicate that as the contrast increases, the general rule becomes progressively less effective at predicting the optimal boundary discretization.

\begin{figure}[!b]
    \centering
    \includegraphics[width=0.8\textwidth]{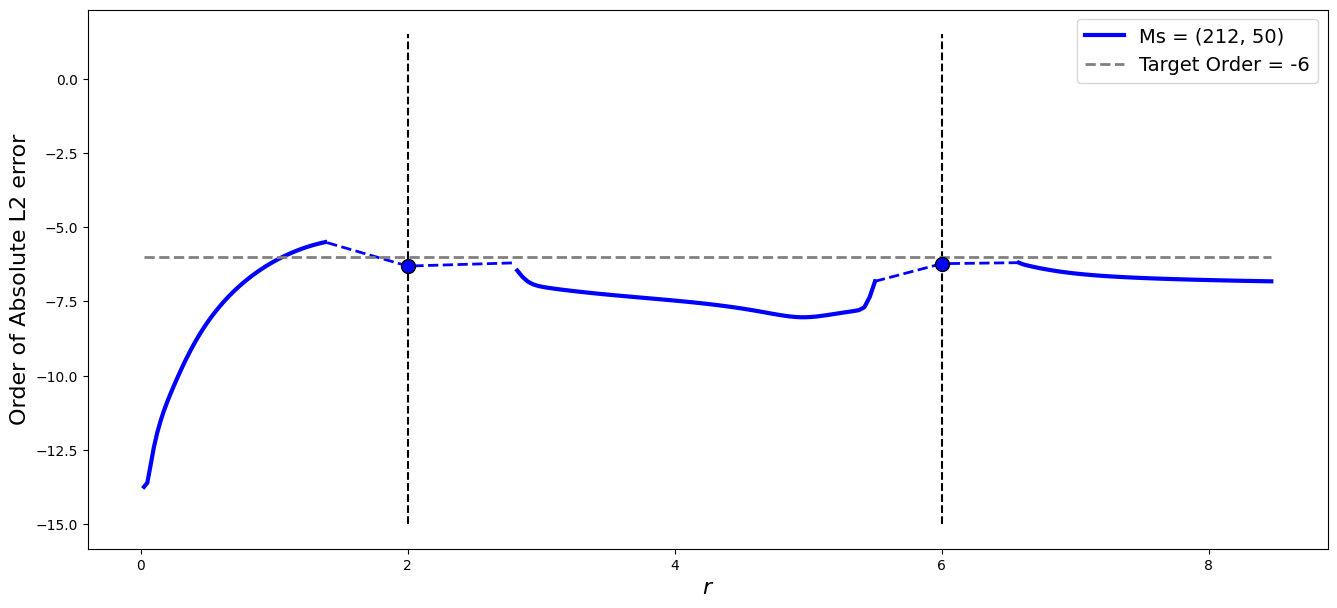}
    \caption{Semilog plot of the absolute \(L_2\) error of BIE solutions case for \(\mathbf k = (8, 2, 6), \mathbf r = (6, 2)\). Only the results from the optimal \(M^{tar}\) values \(\mathbf M = (212, 50)\) are shown.}
    \label{fig:radial errors set 4}
\end{figure}

Choosing an optimal boundary discretization becomes more complex as we introduce additional interfaces \(\Gamma_j\) for the multilayer problem. We present the result for \textbf{Case 5} (with 2 interfaces) in Figure \ref{fig:radial errors set 4}. We obtain two optimal discretizations on \(\Gamma_0\) and \(\Gamma_1\) respectively,  \(\boldsymbol{M^{tar}_0 = 212}\) and \(\boldsymbol{M^{tar}_1 = 50}\), to achieve the chosen accuracy. Although we may attribute the difference in number of points to the varying boundary sizes for \(R_0 > R_1\), the ratio \(\frac{R_0}{R_1} = 3\) does not fully account for the fact that $M^{tar}_0$ is over four times $M^{tar}_1$.  With this boundary discretization, the $L^2$ error seems uniform across the domain.

These numerical experiments suggest that as long as the interfaces are discretized with a number of points large enough to reach the target accuracy, the solutions within each layer can be captured with comparable or better precision. However, this naive adaptation of the discretization has the prohibitive need for the analytic solutions. Furthermore, identifying a rule to account for multiple contrasts \(\beta_j\) and their influence on the interface where they are defined is challenging.  It is even more difficult to quantify their effect on other interfaces. We would like to develop an \textit{automatic adaptive} approach, i.e., a method that is applicable to cases where we do not have access to analytic solutions.

In the following section, we introduce an adaptive mesh strategy that optimizes the discretization to minimize the interpolation error.  Unlike the boundary-fitted approach it refines the mesh only where it is needed, enabling accurate and efficient computation of the field solutions.

\section{Adaptive Mesh  Refinement}
\label{sec:adaptive mesh generation}
\subsection{From Interpolation Error for Layered Domains to Optimal Mesh Design}
\paragraph*{Deriving the optimal volume mesh} We use an automatic mesh adaptation algorithm to determine the optimal discretization for each domain (\textit{Adaptive Method}). 
Our goal is to determine a mesh (with   size $K$) that minimizes the interpolation error, estimated from a truncated Taylor expansion together with geometric considerations. However, solving this minimization problem directly at the discrete level is challenging, which motivates the use of an equivalent continuous formulation. 

We follow the approach of \cite{loseille2011,loseille2017}. An anisotropic adapted mesh $\boldsymbol{T}^*$ is generated with respect to a metric $\mathcal{T}$, so that all elements have unit size in the $\mathcal{T}$-metric. The optimal metric $\mathcal{T}$ is obtained by minimizing the continuous interpolation error for a given mesh complexity $\mathcal{N}$, thereby defining the corresponding adapted mesh $\boldsymbol{T}^*$ and controlling the overall mesh density.

 Given a domain $\Omega$ with smooth closed boundary $\Gamma$, the optimal metric \(\mathcal{T}\), in the \(L^p\) norm (for \(p \geq 1\)), is given by

\begin{equation}
    \mathcal{T} =  \mathcal{N}^\frac{2}{3} \left( \int_\Gamma \det(|H_R(u)|)^{\frac{p}{2p + 3}} \right)^{-\frac{2}{3}} \det(|H_R(u)|)^{\frac{-1}{2p + 3}}|H_R(u)|
\end{equation}
where $H_R(u)$ is the recovered Hessian of the solution.  
In addition, with this approach, an explicit bound on the interpolation error associated with $\boldsymbol{T}^*$ is available and given by

\begin{equation}
    E_p(\boldsymbol{T}^*) = \frac{\mathcal{N}^{-1}}{4} \left( \int_\Gamma \det(|H_R(u)|)^{\frac{p}{2(p+1)}}\right)^{\frac{p+1}{p}}.
    \label{eq:interp error bound}
\end{equation}

Finally, the corresponding error estimate becomes

\begin{equation*}
  \frac{1}{2} E_p(\boldsymbol{T}^*) \leq  || u - \Pi_h u||_{L^p(\Omega)} \leq 2E_p(\boldsymbol{T}^*)
\end{equation*}
where  \(\Pi_h u\) is  the linear interpolant of the solution \(u\) on the mesh $\boldsymbol{T}^*$.

\paragraph*{Illustration of the behavior of the Adaptive Method}

For each interface $\Gamma_j$, we begin by defining a fine discretization: a dense\footnote{In general the algorithm may start from a coarse mesh, but in this section a fine mesh is used for preserving the interfaces' geometry.}, equispaced set of candidate points on the boundary.  These points serve as part of the initial data for our mesh adaptation process.  From there, we refine or coarsen the mesh based on a mesh complexity input requirement.  However,  a key constraint is imposed: the interface geometry must be preserved.  This ensures the interface remains geometrically consistent and is not altered during mesh refinement.

Although we begin with an equispaced set of points on the interfaces, the algorithm does not necessarily preserve that spacing.  Thus, the distribution of retained interface vertices may be nonuniform.  

\begin{figure}[!h]
\centering
\begin{subfigure}[t]{0.32\textwidth}
        \centering
        \includegraphics[width=0.95\textwidth]{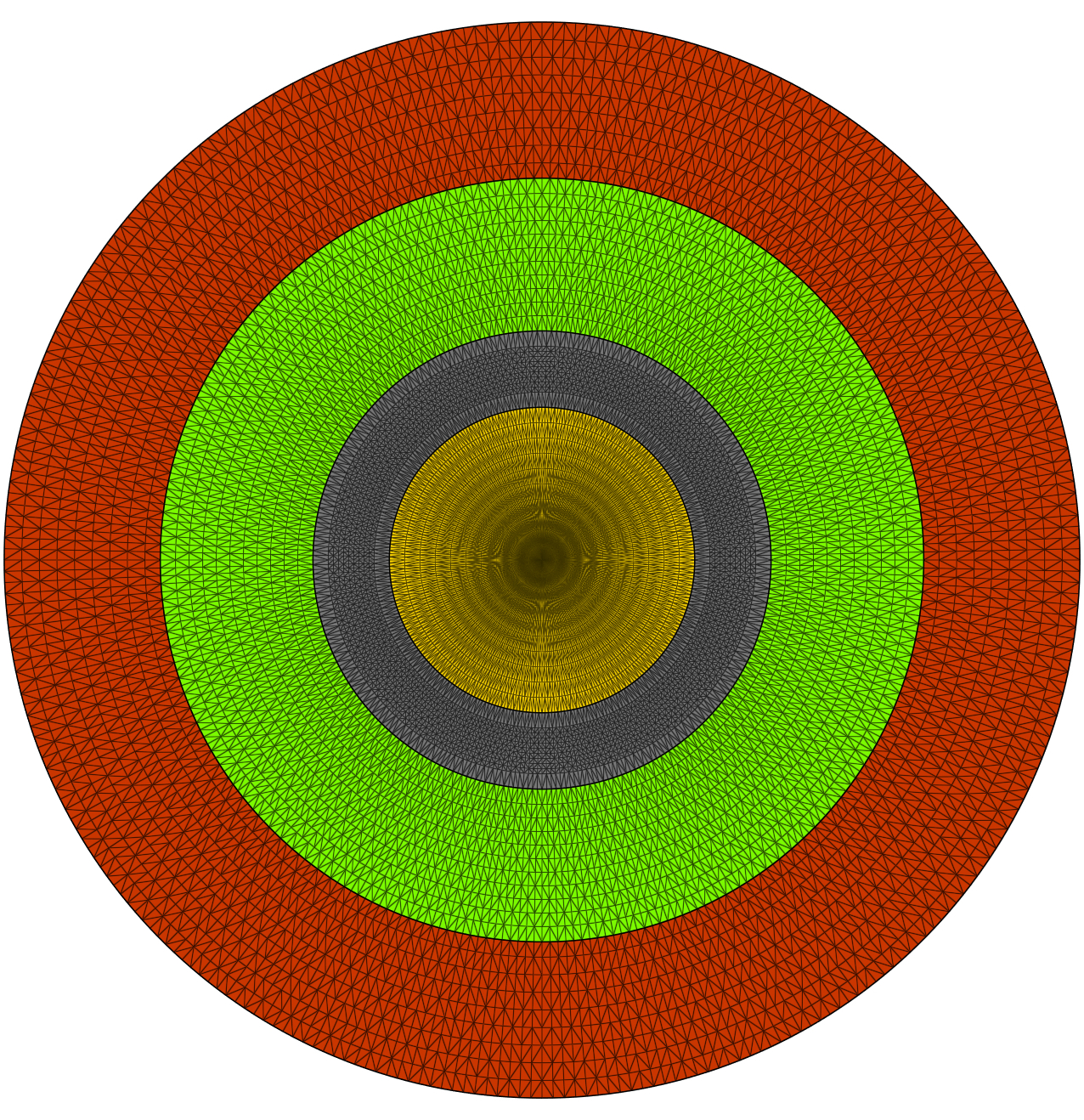}
        \caption*{(a) Initial, boundary-fitted discretization.}
\end{subfigure}\hspace{1mm}
\begin{subfigure}[t]{0.32\textwidth}
        \centering
        \includegraphics[width=0.95\textwidth]{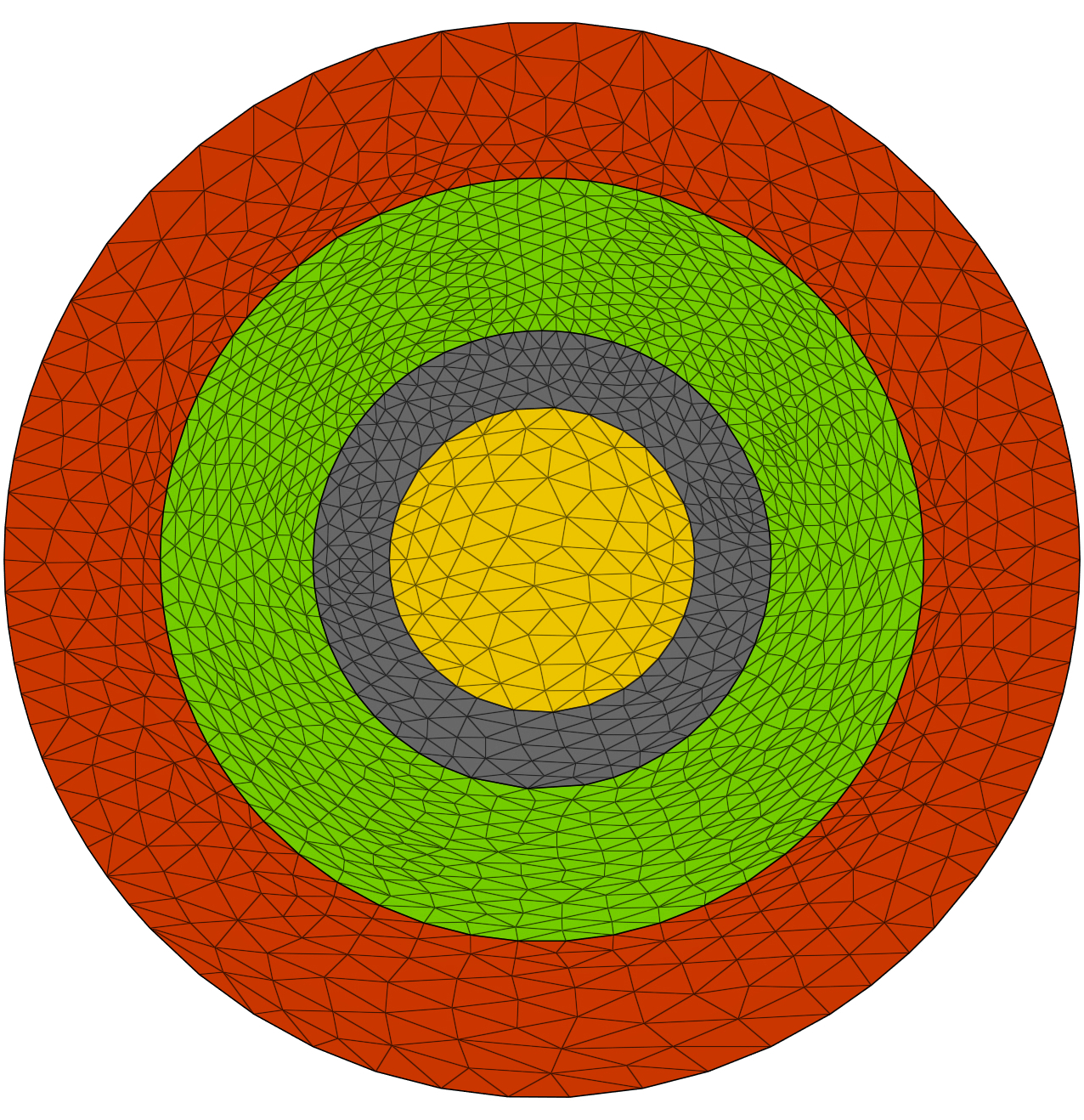}
        \caption*{(b) Adaptive Mesh.}
\end{subfigure}\hspace{1mm}
\begin{subfigure}[t]{0.32\textwidth}
        \centering
        \includegraphics[width=0.95\textwidth]{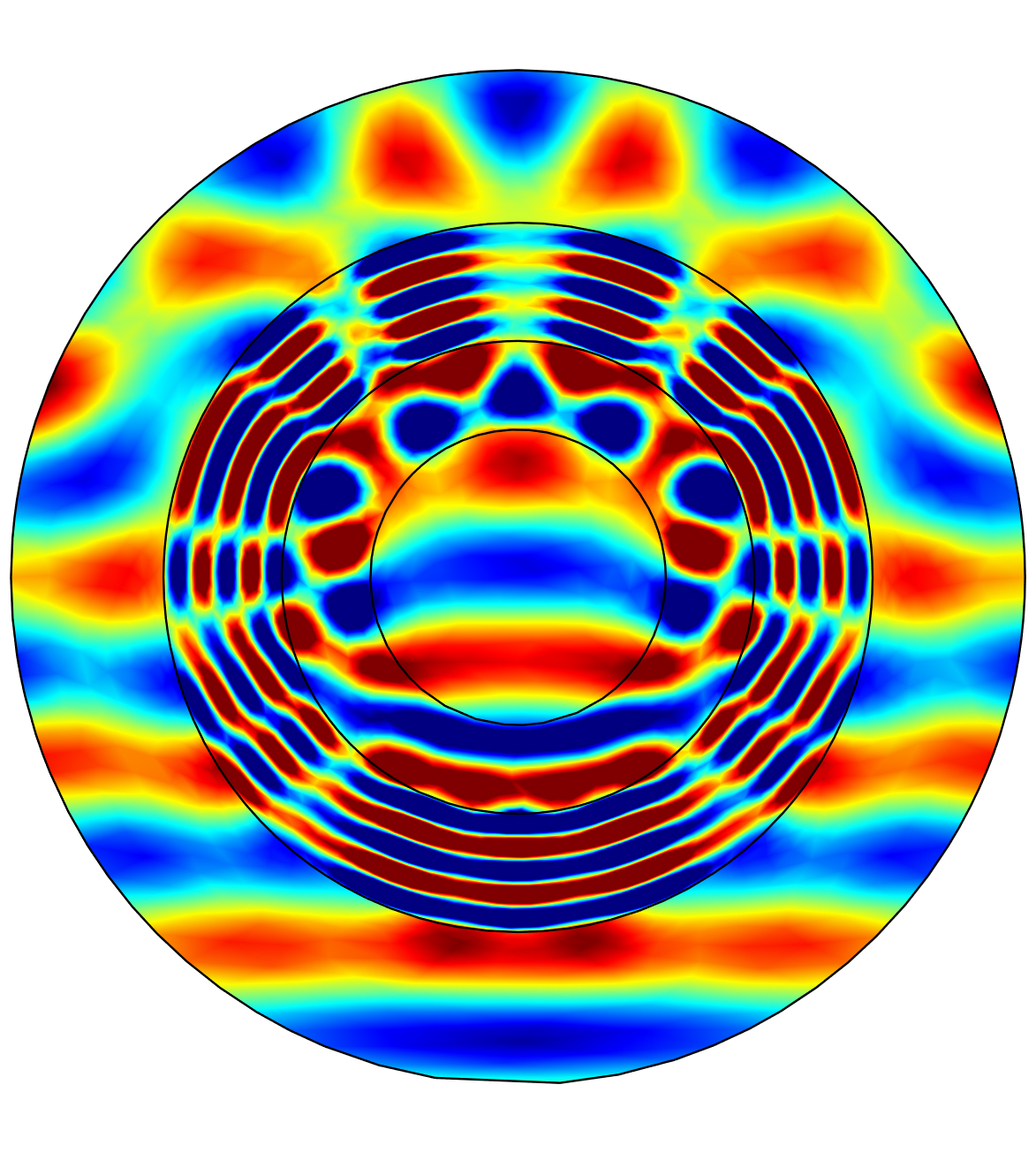}
        \caption*{(c) Real part of the solution.}
\end{subfigure}
        \caption{Comparison of meshes at initial and later iterations of mesh generation versus the solution for \(\mathbf{r} = (6,4,2.5)\) and \(\mathbf{k} = (2,8,4,2)\).}
    \label{fig:mesh_comparisons}
    \end{figure}

As an illustration, Figure~\ref{fig:mesh_comparisons} presents a representative comparison between a boundary-fitted mesh and its adapted (anisotropic) counterpart for a configuration with three interfaces and varying material contrasts. As expected, the adaptive mesh concentrates points in regions where higher resolution is required, leading to highly nonuniform discretizations along each interface compared with the uniform spacing observed in Figure~\ref{fig:mesh_comparisons}(a).

\paragraph*{Automatic mesh adaptation}

The adaptive method follows an iterative process: we initialize the algorithm with a coarse target complexity, typically \(\mathcal{N}_0 = 200\), and compute the solution in the volume.  From this solution, we estimate the Hessian and construct the corresponding optimal metric \(\mathcal{T}\) as well as the corresponding estimate of the interpolation error $E_p(\boldsymbol{T}^*)$.  A new mesh is then generated from this metric for a chosen \(\mathcal{N}\), and this process is repeated.

Rather than increasing the mesh complexity at every iteration, we raise it only intermittently, e.g., by a fixed percentage every other iteration.  This strategy offers the possibility to optimize the interpolation error at a given complexity before refining the mesh.  At each iteration, the interpolation error associated with the current metric is measured.  The relative improvement between iterations serves as an indicator of whether further adaptation at the same complexity is needed.  The algorithm may be stopped once the improvement in the interpolation error becomes negligible or when a prescribed error bound estimate or maximum complexity is reached.

The adaptive method is coupled with a two-stage procedure to evaluate solutions.  We first find the boundary data by solving the system of corresponding boundary integral equations \eqref{eq:boundaryreps}.  The solution is then propagated into the volume via boundary integral representation  \eqref{eq:bie_layers}.  To utilize boundary element methods, the boundary geometry is preserved by the adaptive method.

A particular constraint of our implementation arises from our choice of quadrature rule used for the boundary integrals: they require a discretization that is equi-spaced, defined at \(M_j\) values of \(\theta_m \) \eqref{eq:theta distribution}.  To satisfy this requirement, at each iteration we extract from the current mesh its closest approximation to a uniform boundary mesh.  This is used to evaluate accurate boundary data for the solutions.  The volume mesh, however, is adapted strictly according to the optimal metric and is not required to be uniform.  It is our chosen quadrature that requires an equi-spaced distribution.  In return it yields spectral accuracy in the volume, so as long as boundary solutions are sufficiently accurate, the solution remains highly accurate throughout the domain (as shown in section 3.3).  This enables accurate wave field solutions with reduced computational cost for reduced, optimal discretizations.

\subsection{Illustration of the adaptive process on circular boundaries}

We now regard the mesh adaptation process on a series of representative examples involving circular boundaries (presented in Table~\ref{tab:multbdycases}). Importantly, we want to determine whether a general rule such as \eqref{eq:general rule} is appropriate.

\begin{table}[!h]
\centering
\begin{tabular}{|c c c c|}
\hline
Case & $\mathbf{k} = (k_0, k_1,...,k_N)$ & $\mathbf{r} = (R_0, ..., R_{N-1})$& $\mathbf{b} = (\beta_0, ..., \beta_{N-1})$\\
\hline
1 & $(2,6)$  & $4$ & 9\\ 
5& $(2,6,1)$& $(4,2)$& (9, $\tfrac{1}{36}$)\\ 
6& $(2,6,10)$& $(4,2)$& $(9, \tfrac{25}{9})$\\
7& $(6, 2, 4, 1)$& $(4,2,1)$& $(\tfrac{1}{9}, 4, \tfrac{1}{16})$\\
\hline
\end{tabular}
\caption{Summary of multiple-boundary cases used in numerical experiments with the Adaptive Method, with wavenumbers $\mathbf{k}$, boundary radii \(\mathbf{r}\), and contrasts \(\mathbf{b}\).}
\label{tab:multbdycases}
\end{table}

We first consider \textbf{Case 1}, for which  we already gained some knowledge in previous section.   

\begin{figure}[!b]
\centering
        \includegraphics[width=\linewidth]{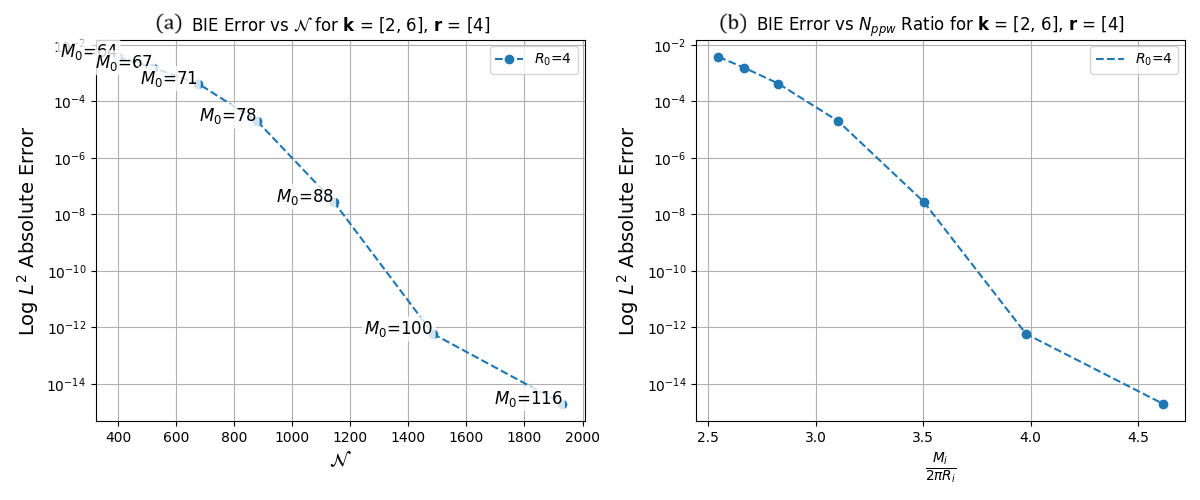}
        \caption{For \textbf{Case 1}: \(\mathbf k = (2, 6)\) and \(R_0 = 4\): (a) Log of the absolute \(L_2\) error \eqref{eq:error_approx} obtained with the BIE solution on \(\Gamma_0\) versus the mesh complexity \(\mathcal{N}\), for various boundary resolutions \(M_0\) as labeled. (b) The same error plotted against \(\frac{M_0}{2\pi R_0}\).}
        \label{fig:1boundaryerrors}
    \end{figure}

Recall that the adaptive method relies on an iterative process. To obtain the results presented in Figure \ref{fig:1boundaryerrors}, we fix the initial complexity across the complete domain at \(\mathcal{N}_0 = 400\).   The automatic mesh adaptation process undergoes 14 iterations, where we increase the complexity $\mathcal{N}_n$ by 30\% every other iteration and otherwise re-mesh on the same complexity once.  As the mesh refines, the number of boundary points \(M_0\) changes.  

We have an adapted mesh, but we look at the accuracy of the BIE solutions by computing the corresponding absolute \(L^2\) error using the known analytic solution as in \eqref{eq:error_approx} at each iteration.   Figure~\ref{fig:1boundaryerrors} shows the evolution of this error as a function of the mesh complexity $\mathcal{N}$ and the number of boundary points $M_0$. 

Examining Figure~\ref{fig:1boundaryerrors}(a), we observe that the \(L^2\) error is highly sensitive to the choice of \(M_0\).  In particular, looking at where the discretization is refined by approximately 10 additional points, increasing \(M_0\) from 78 to 88 decreases the error by roughly one order of magnitude, while increasing \(M_0\) to \(98\approx100\) results in a reduction of several orders.  This highlights that even slight deviations from the optimal number of points, particularly under-resolving, can have substantial impacts on the BIE solution's accuracy. 

\begin{figure}[!b]
    \centering
    \begin{subfigure}{\textwidth}
        \centering
        \includegraphics[width=0.95\textwidth]{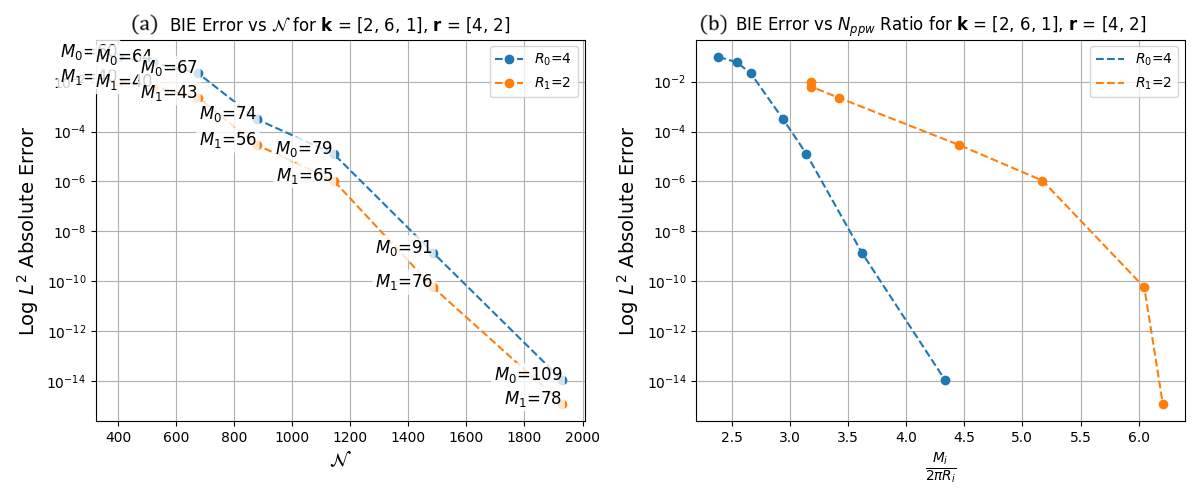}
        \caption*{\textbf{Case 5}: \(\mathbf k = (2, 6, 1), \mathbf r = (4, 2)\)}
    \end{subfigure} 
    \begin{subfigure}{\textwidth}
        \centering
        \includegraphics[width=0.95\textwidth]{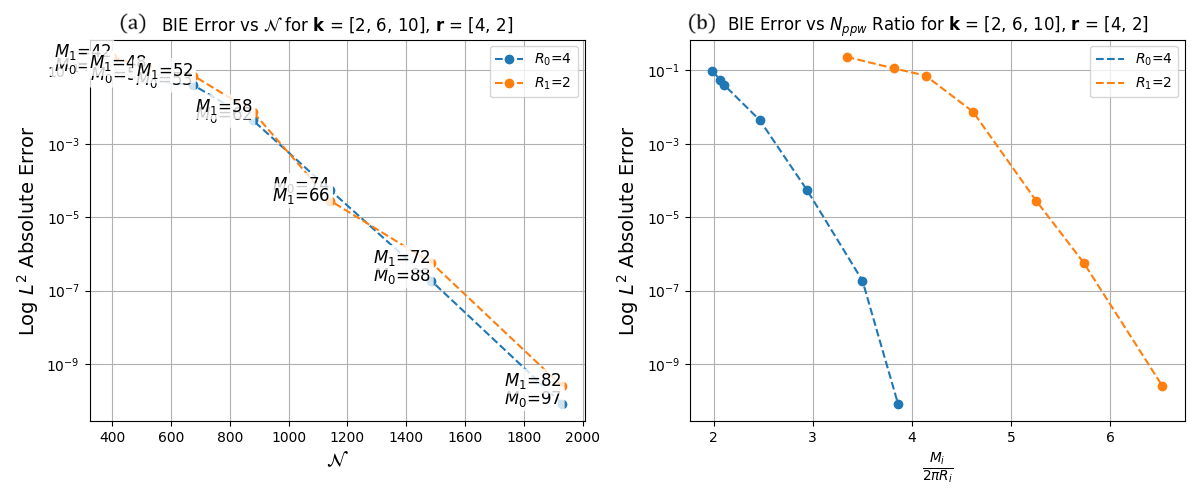}
        \caption*{\textbf{Case 6}: \(\mathbf k = (2, 6, 10), \mathbf r = (4, 2)\)}
    \end{subfigure}
    
   \caption{(a) Log absolute \(L_2\) error of the BIE solutions on \(\Gamma_0\) and \(\Gamma_1\)versus total mesh complexity \(N\), for various boundary resolutions \(M_0\) and $M_1$ as labeled. (b) The same error plotted against the ratio \(\tfrac{M_i}{2\pi R_i}\).}
        \label{fig:2boundaryerrors}
\end{figure}
\begin{figure}[!b]
\centering
        \includegraphics[width=\linewidth]{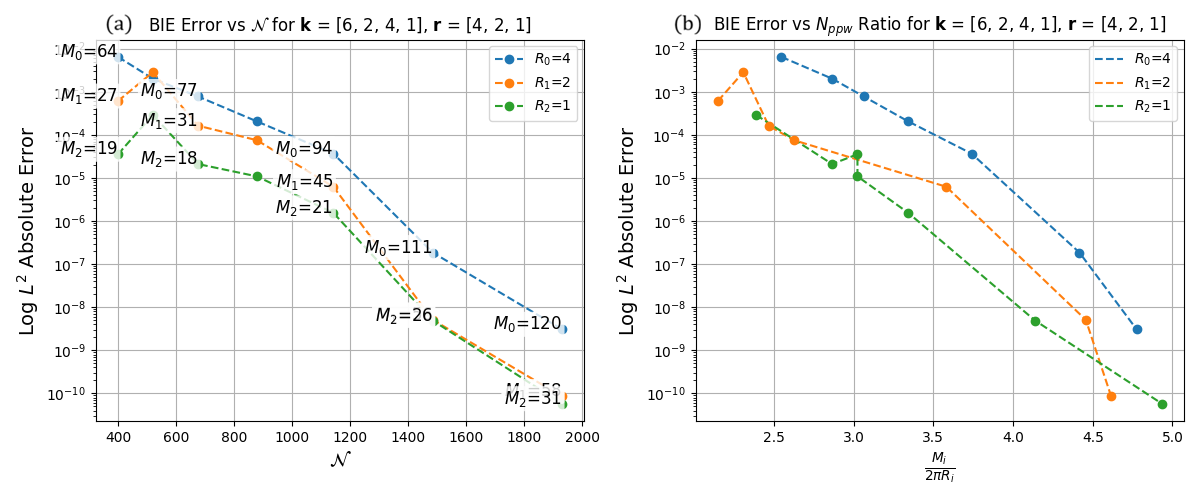}
        \caption{For \textbf{Case 7}: \(\mathbf k = (6, 2, 4, 1)\) and \(\mathbf{r} = (4, 2, 1)\): (a) Log absolute \(L_2\) error of the BIE solution on three boundaries versus total mesh complexity \(N\), for various boundary resolutions \(M_0\), \(M_1\), and \(M_2\) as labeled. (b) The same errors plotted against the ratio \(\tfrac{M_i}{2\pi R_i}\).}
        \label{fig:3boundaryerrors}
\end{figure}

We recall that our goal is to verify the accuracy of the general rule. 
Figure~\ref{fig:1boundaryerrors}(b) plots the error versus the ratio $\frac{M_0}{2\pi R_0}$, which is proportional to the points-per-wavelength metric $N_{\rm ppw}$ from \eqref{eq:general rule}. The results show that this ratio is not constant; if it were, it would provide a simple guideline for selecting the optimal boundary discretization.  To achieve a target accuracy of $10^{-6}$, the adaptive mesh algorithm reaches an optimal resolution with only $\frac{M_0}{2\pi R_0} \approx 3$, well below the typical rule of 5–10 points per wavelength for $k_{max} = 6$.  This can be explained by the lower exterior wavenumber $k_0 = 2$, which dominates the mesh requirements and reduces the number of boundary points needed for $k_1 = 6$. Even for errors on the order of $10^{-14}$, the ratio remains below the standard 5–10 points-per-wavelength guideline.

After this introductory example, we turn to more complex cases with additional circular boundaries, as given in Table \ref{tab:multbdycases}.

Figure \ref{fig:2boundaryerrors} presents results for the 2-boundary cases from Table \ref{tab:multbdycases} (\textbf{Case 5} and \textbf{Case 6}). The Adaptive Method achieves a comparable order of accuracy on both boundaries, $\Gamma_0$ and $\Gamma_1$, with different discretizations. Again, in both cases there is no clear trend when comparing to the ratios $\tfrac{M_i}{2\pi R_i}$ ($i = 0, 1$), reinforcing the absence of a general rule for optimal boundary discretization in multilayer transmission problems. Increasing $k_{max}$ from 6 to 10 in \textbf{Case 6} results in only slightly increased discretization demand on $\Gamma_1$ while it remains relatively unchanged on $\Gamma_0$. The discretization range that might otherwise be assumed from the general rule \eqref{eq:general rule} for \(5 \leq N_{ppw} \leq 10\), given \(k_{max} = 10\), is between \(M = 200\) to \(400\).

Figure~\ref{fig:3boundaryerrors} shows that similar trends are observed for \textbf{Case 7}, with three boundaries.  Additionally, though a maximum wavenumber of \(k_0 = 6\) is introduced in the outermost subdomain, the presence of subdomains with lower wavenumbers substantially reduces the overall discretization requirements on the interior interfaces.  The transmission conditions couple the interface behavior in such a way that improved accuracy on \(\Gamma_1\) and \(\Gamma_2\) is obtained with little additional discretization effort, except on \(\Gamma_0\).  Notably, increasing the discretization on \(\Gamma_0\) to improve its boundary accuracy also improves the accuracy on the other boundaries.

From these circular boundary examples examined, the Adaptive Method seems to maintain similar accuracy across all boundaries by concentrating resolution where needed, often reducing mesh requirements due to lower wavenumbers in sub-regions $\Omega_j$. Although no simple rule applies universally, it reliably identifies boundary discretizations that preserve uniform accuracy, even in multilayer configurations with high material contrasts.

Having validated the Adaptive Method on simple geometries where an analytical solution can be found, we now turn to more complex boundaries. In the next section, we assess its performance on \textit{exotic} shapes, where analytic solutions are unavailable and local refinement, curvature, and material contrasts pose greater challenges.

\subsection{Adaptive Mesh Algorithm with BIE Solutions}

\begin{table}[!b]
\centering
\begin{tabular}{|c l c|}
\hline
\textbf{Method}& \textbf{Description} & \textbf{Computation}\\
 &\vspace{-0.3cm}\\
\hline
\textbf{ADAPT-ANA}& Interpolation error bound for &\vspace{-0.3cm} $E_2(\boldsymbol{T}^*_{ana})$ \\
 & \vspace{-0.3cm}\\&  analytic solution on adapted mesh $\boldsymbol{T}^*_{ana}$ & \\
 \hline 
\textbf{ADAPT-BIE}& Interpolation error bound for BIE solution & \vspace{-0.3cm} $E_2(\boldsymbol{T}^*)$\\
 &\vspace{-0.3cm}\\
 &on adapted mesh $\boldsymbol{T}^*$ & \\
 \hline
\end{tabular}
\caption{Description of methods plotted in Figure \ref{fig: feflo errors 2}.}
\label{tab:methods}
\end{table}

In general, one does not have access to an analytic solution to investigate how the Adaptive Method is performing. Since the volumetric interpolation error bound \eqref{eq:interp error bound} should consistently improve over increasing mesh complexities if the algorithm is performing well, we simply rely on computing \eqref{eq:interp error bound} using the approximated solution, $u$. Here we make use of the BIE solutions and the boundary integral representation to obtain the recovered Hessian $H_R(u)$.

In order to validate our approach, we first test on circular boundaries so we can compare with previous results. We call \textbf{ADAPT-ANA} the error bound in terms of the recovered Hessian \(H_R(u^{ana})\) which defines the mesh \(\boldsymbol{T}^*_{ana}\), and \textbf{ADAPT-BIE} the one where the recovered Hessian is determined from the BIE solution, defining the mesh \(\boldsymbol{T}^*\). We summarize the two approaches in Table \ref{tab:methods}.

\begin{figure}[!h]
\centering
        \includegraphics[width=0.495\textwidth]{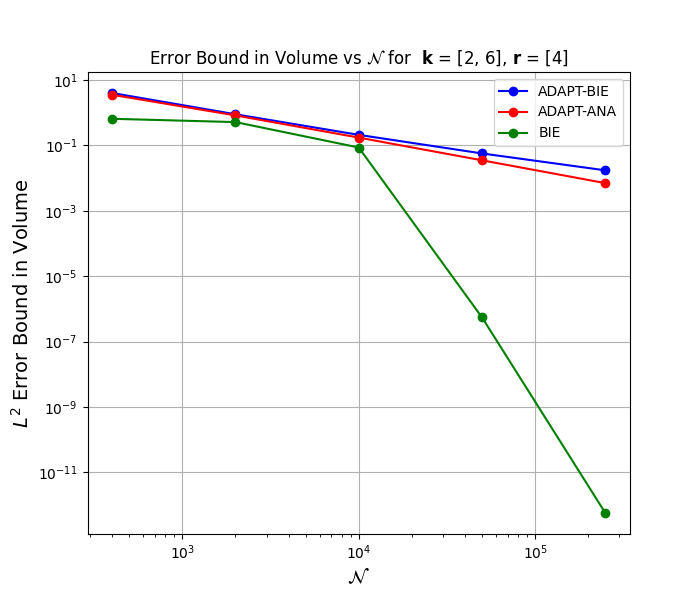} \includegraphics[width=0.495\textwidth]{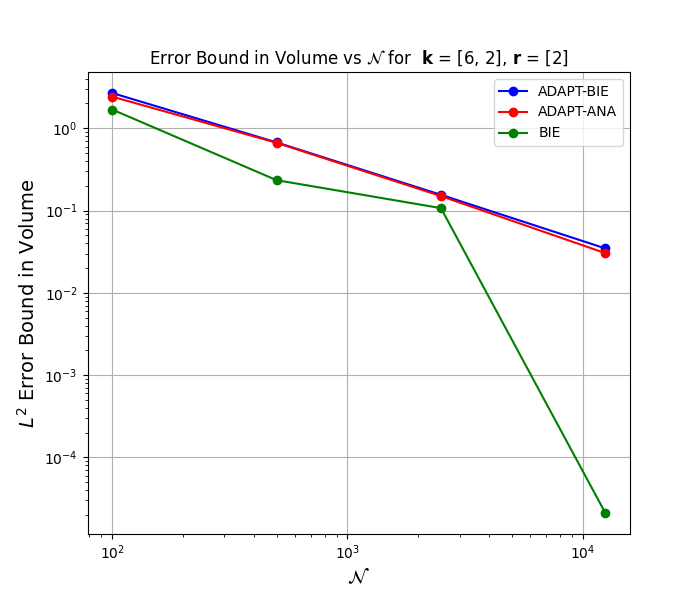}
        \caption{Plot of the error bounds described in Table \ref{tab:methods} for \textbf{Case 1 and 2}, over increasing mesh complexities $\mathcal{N}$.}
        \label{fig: feflo errors 2}
\end{figure}

Figure \ref{fig: feflo errors 2} presents the interpolation error bounds computed with both methods for mesh complexities increasing from \(10^3\) to \(10^5\) for \textbf{Case 1} and \textbf{Case 2}. We find linear convergence with \textbf{ADAPT-ANA}, which is in accordance with theoretical expectations \cite{chaillat_metric-based_2018}. We observe convergence with \textbf{ADAPT-BIE}, however with a slower rate for larger mesh complexities. This behavior is likely caused by increased refinement near the boundaries, where the evaluation of the boundary potentials becomes challenging due to nearly-singular kernels.  This \textit{close-evaluation } error is well-known and can impact the convergence behaviors of standard quadrature rules \cite{barnett_evaluation_2014}.  This highlights an important research direction when combining adaptive meshing with BEM solvers: the integration of close-evaluation treatment methods.

The first two error estimates give an upper bound but we would like to check if we can rely on it.  Here, we also estimate the error in the volume by comparing the BIE solution with the analytic solution: this is what we call \textbf{BIE} in Figure \ref{fig: feflo errors 2}. This error is computed via the $L^2$ error 
\begin{equation}
    ||u - u^{ana}||_{L^2(\tilde{\Omega})} =\sum \limits_{j} ||u_j - u_j^{ana}||_{L^2(\Omega_j \setminus {\omega_j})} 
\end{equation}
where $\tilde{\Omega}$ is the ``perforated" domain $\Omega$, where we remove interfaces vicinities (called $\omega_j$) to avoid the so-called close evaluation problem. This yields a more accurate measure of the accuracy of the solution in the volume as we do not treat the nearly-singular behaviors. Figure \ref{fig: feflo errors 2} shows spectral convergence of \textbf{BIE}, as expected theoretically. 

\begin{figure}[!b]
    \centering
    \begin{subfigure}{\textwidth}
        \centering
        \includegraphics[height=6cm,width=0.8\textwidth]{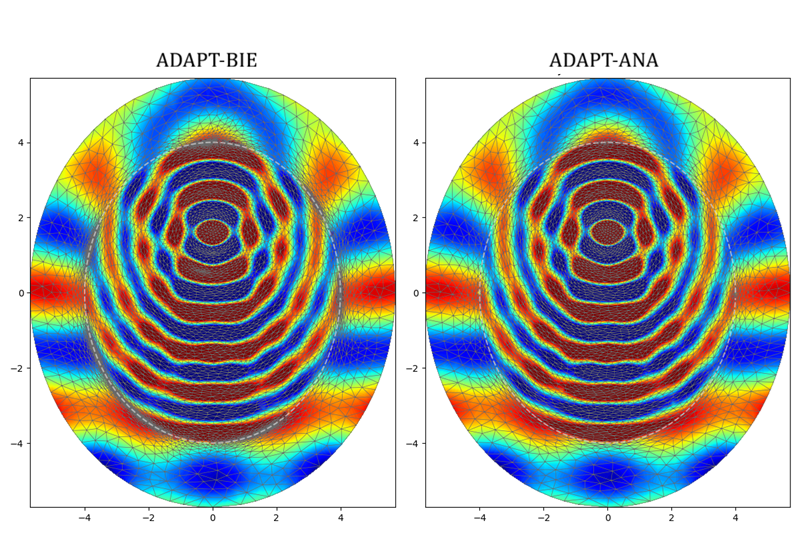}
        \caption*{(a) \(\mathcal{N} = 2 \cdot 10^3\) }
    \end{subfigure} 
    \begin{subfigure}{\textwidth}
        \centering
        \includegraphics[height=6cm,width=0.8\textwidth]{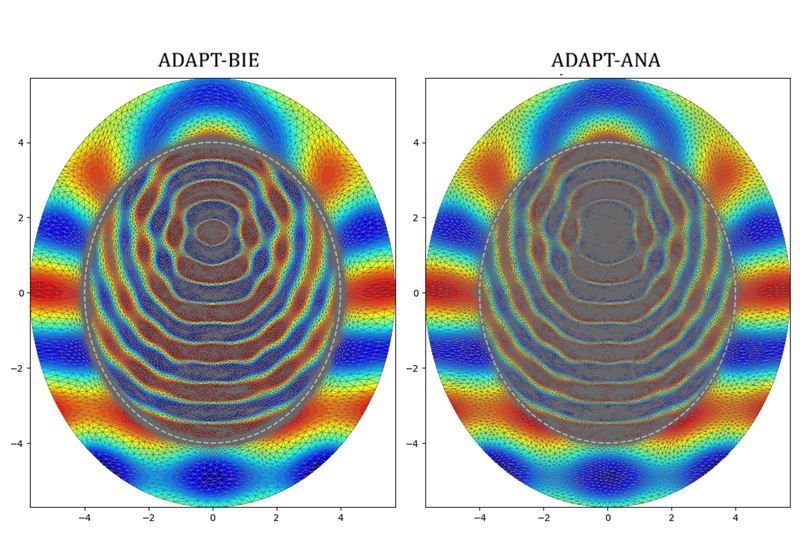}
        \caption*{(b) \(\mathcal{N} = 3.2 \cdot 10^4\)}
    \end{subfigure}
   \caption{Comparison of ADAPT-BIE and ADAPT-ANA meshes for varying mesh complexities  \(\mathcal{N}\), with \(\mathbf{k} = (2, 6)\) and \(R_0 = 4\) after (a) 5 iterations versus (b) 11 iterations of adaptation.} \label{fig:ptr vs ana meshes}
\end{figure}

For \textbf{Case 1} in particular, we see that the error in the volume begins to decrease exponentially at higher orders of complexity \(\mathcal{N}\). This can likely be attributed to the increasing number of nodes discarded due to their proximity to the boundaries.  At very high complexities, the Adaptive Method places additional nodes near the boundaries to resolve the jumps in the solution magnitude caused by the nearly-singular behavior in those regions.  It is more prominent for the case when \(k_0 = 2\) because the lower exterior wavenumber drives lower discretization demand to capture the overall solution behavior in the domain, so the adaptive method focuses on close evaluation behavior near the boundary.  Figure \ref{fig:ptr vs ana meshes} compares optimal meshes (according to the mesh adaptation algorithm) from ADAPT-BIE versus ADAPT-ANA.  In Figure \ref{fig:ptr vs ana meshes}(a), the mesh densities are similar and relatively uniform through each region \(\Omega_0\) and \(\Omega_1\) in the volume. In Figure \ref{fig:ptr vs ana meshes}(b), for a higher complexity, the ADAPT-BIE mesh shows increased node concentration near the boundary, while the ADAPT-ANA mesh adds more nodes throughout the whole interior region \(\Omega_1\). This near-boundary refinement difference only becomes significant at sufficiently high complexities relative to the domain size.  In this example, it does not appear until \(\mathcal{N} > 10^4\).  Therefore, our results should remain optimal when using BIE solutions for the cases examined so far, and subsequent cases, provided the mesh complexity stays below this threshold and the relative domain size is not reduced.

Having established the interpolation error bound may tell us the effectiveness of the adaptive process where we do not have reference solutions, we now turn to more complex geometries to assess the robustness of the mesh adaptation strategy. 

We use the following parameterization to model star shape objects:
\begin{equation}
\label{eq:star_shape}
\begin{aligned}
x_i(\theta) &= \big(R_i + a \cos(n \theta)\big)\cos \theta, \\
y_i(\theta) &= \big(R_i + a \cos(n \theta)\big)\sin \theta,
\end{aligned}
\qquad i = 0, \ldots, N-1,
\end{equation}
where \(n\) is the number of star points and \(a \geq 0\) controls their relative size with respect to \(R_i\). 

We can leverage the error bound associated with the BIE solutions to determine whether additional re-meshing iterations are necessary, and how many. If the error bound continues to improve with each iteration, the current number of re-meshing iterations is likely sufficient.  Conversely, if the error bound plateaus or degrades, this signals the need to adjust the number of re-meshing iterations to achieve better results.

\begin{figure}[!h]
    \centering
    \hfill
    \begin{subfigure}[t]{0.45\textwidth}
        \centering
        \includegraphics[width=\textwidth]{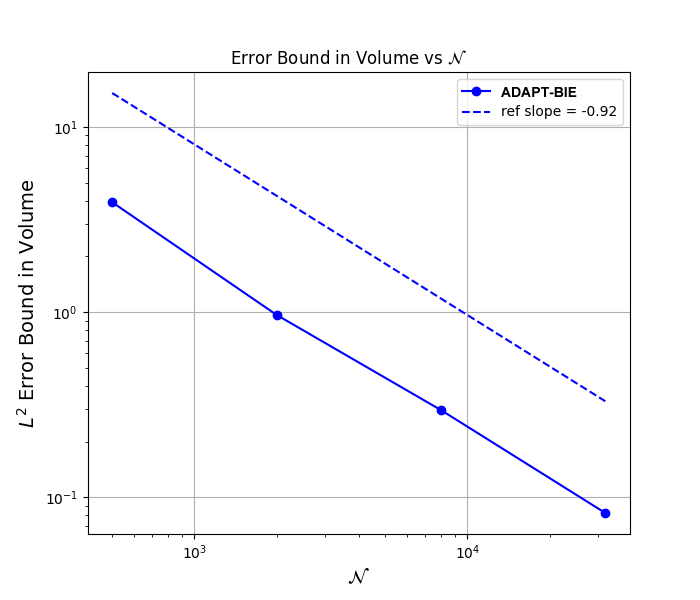}
        \caption*{(a) No re-meshing before refinement}
    \end{subfigure}
    \begin{subfigure}[t]{0.45\textwidth}
        \centering
        \includegraphics[width=\textwidth]{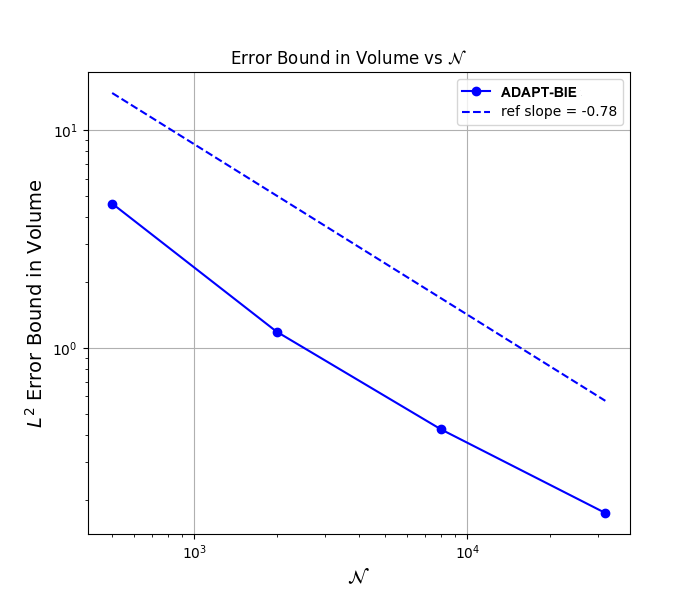}
        \caption*{(b) Two re-meshings before refinement}
    \end{subfigure}
    \caption{Error bounds and reference slopes of the interpolation versus complexity \(\mathcal{N}\) for \(\mathbf{k} = (2, 6)\), \(R_0 = 4\), \(a = 0.05\), \(n = 10\) for different re-meshing schemes.}
    \label{fig:remeshing error bound a=0.05}
\end{figure}

\begin{figure}[!h]
    \centering
    \hfill
    \begin{subfigure}[t]{0.45\textwidth}
        \centering
        \includegraphics[width=\textwidth]{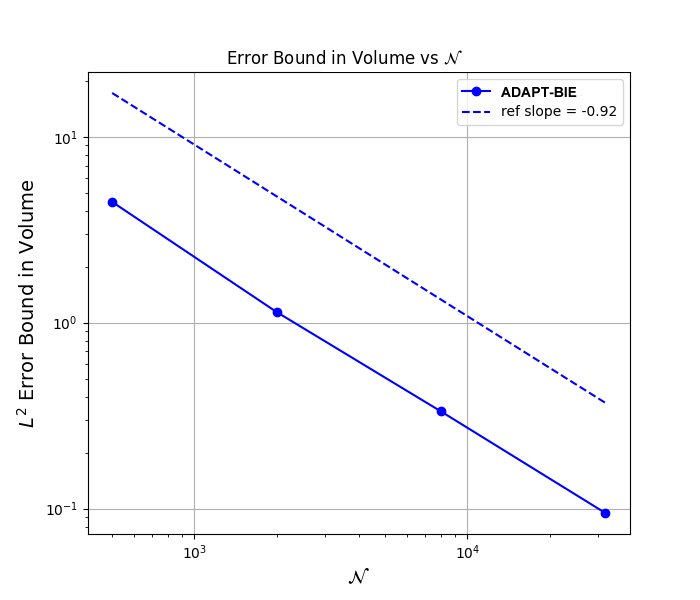}
        \caption*{(a) No re-meshing before refinement}
    \end{subfigure}
    \begin{subfigure}[t]{0.45\textwidth}
        \centering
        \includegraphics[width=\textwidth]{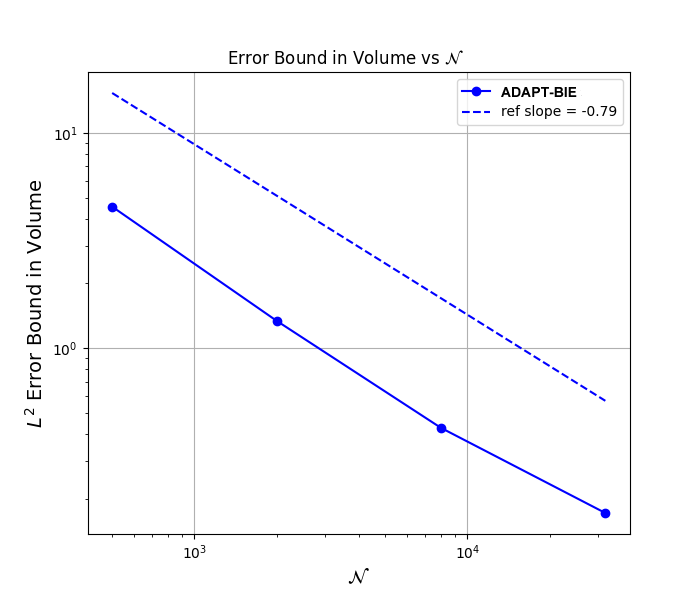}
        \caption*{(b) Two re-meshings before refinement}
    \end{subfigure}
    \caption{Error bounds and reference slopes of the interpolation versus complexity \(\mathcal{N}\) for \(\mathbf{k} = (2, 6)\), \(R_0 = 4\), \(a = 0.1\), \(n = 10\) for different re-meshing schemes.}
    \label{fig:remeshing error bound a=0.1}
\end{figure}

We examine two cases for amplitudes \(a = 0.05\) and \(a = 0.1\) with n = 10, to compare the behavior of a slight versus more pronounced star point amplitude under two meshing strategies: no re-meshing versus two re-meshings, as shown in Figures \ref{fig:remeshing error bound a=0.05} and \ref{fig:remeshing error bound a=0.1}.  We observe that fewer re-meshing iterations tend to yield better convergence of the error bound. 

\begin{figure}[!h]
    \centering
    \hfill
    \begin{subfigure}[t]{0.49\textwidth}
        \centering
        \includegraphics[width=\textwidth]{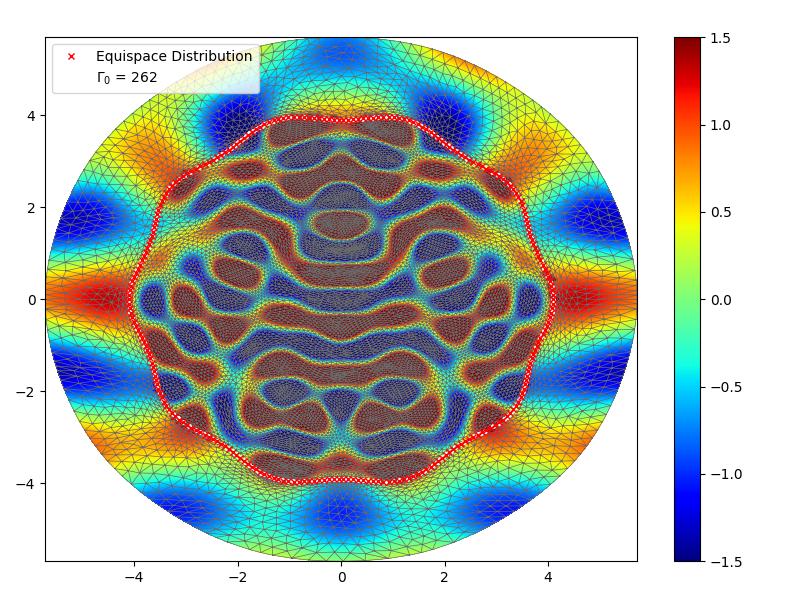}
        \caption*{(a) Initial mesh generated for \(\mathcal{N} = 2000\)}
    \end{subfigure}
    \begin{subfigure}[t]{0.5\textwidth}
        \centering
        \includegraphics[width=\textwidth]{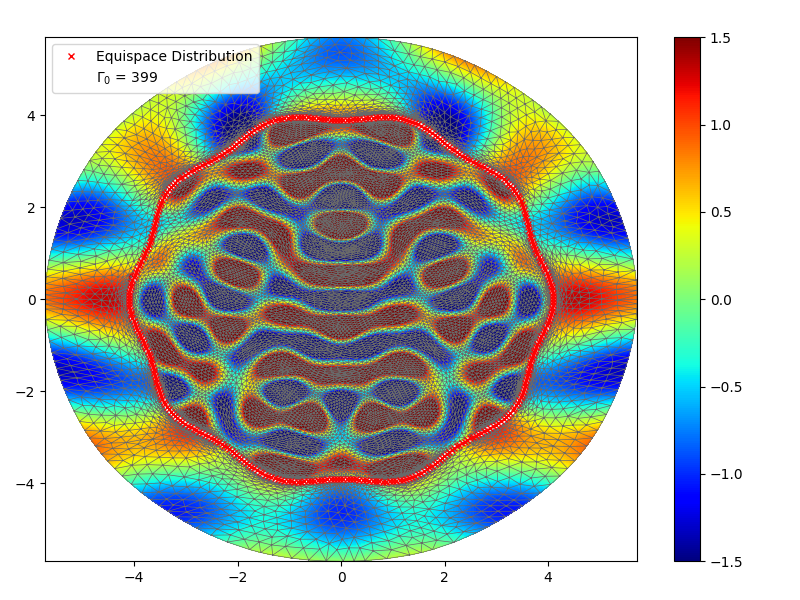}
        \caption*{(b) Mesh generated after re-meshing twice for \(\mathcal{N} = 2000\)}
    \end{subfigure}
    \caption{Real part of the PTR solutions and generated meshes for \(\mathbf{k} = (2, 6)\), \(R_0 = 4\), \(a = 0.1\), \(n = 10\).  We note for the same mesh complexity \(\mathcal{N} =  2000\), after re-meshing from (a) the new mesh generated in (b) adds more refinement to the boundary.}
    \label{fig:mesh vs remeshed stars}
\end{figure}

For these cases, we also consider that the adaptive strategy may prioritize regions of high curvature and sharper geometries, allocating more points in these regions where interpolation error is likely to be higher regardless of the contrasts.  In particular, cases with high exterior wavenumbers consistently demand finer resolution across all boundaries, while otherwise lower exterior wavenumbers lessen these demands.  Since we cannot quantify the error on the boundaries directly in the absence of analytic solutions, we rely on visual inspection of the corresponding BIE solutions to comment on whether the meshes produced are sufficient to capture the key solution features.  Typically they do, though there is improvements to be made to leverage the same optimality of results we achieve from analytic solutions.  Treatment of close evaluation error near the boundaries while processing solutions, for example, will improve usage of this adaptive algorithm by improving the accuracy of the method. It will also allow re-meshing to be used as intended, to optimally redistribute the mesh points to capture the actual solution behavior rather than try to resolve problematic solution features.  With treated solutions, the adaptive method is expected to limit unnecessary over-resolution near the boundaries.  Although some refinement is unavoidable due to boundary curvature in the cases shown in this section, untreated solutions likely amplify this effect, as shown in Figure \ref{fig:mesh vs remeshed stars}. 

Our qualitative observations reinforce the key insight from the previous section: there is no simple rule linking the optimal number of points $M_i$ on  $\Gamma_i$ to the maximum wavenumber.  Instead, the discretization depends on both the local geometry and material properties within the overall multilayered structure.  However, as the boundary shapes become more complex, the algorithm continues to produce stable and visually reasonable discretizations.  In particular, the background medium, \(k_0\),  strongly influences the discretization behavior.  This complexity highlights the value of automatic boundary discretization approaches in multilayer transmission problems.

\section{Conclusions and Perspectives}

This work has shown that achieving optimal interface discretizations in multilayer transmission problems requires mesh adaptivity, as no simple linear rule—whether based on the contrasts $\beta_j$ or the highest wavenumber—can reliably determine optimal discretizations. Such adaptivity proves essential to obtain the best cost–accuracy balance, particularly due to the interaction between high- and low-wavenumber regimes across material interfaces.

In particular, we have observed that the exterior wavenumber has a significant influence on the overall solution behavior.  A very low exterior wavenumber can damp oscillations within the interior media, whereas a high wavenumber requires a correspondingly finer resolution \cite{liu_bem_2019}.  Because the interfaces often differ not only in physical behavior, size, and shape, using the same number of boundary elements or scaling all interface meshes in proportion to the boundary with the highest wavenumber is not ideal.  Such uniform or proportional discretizations increase computational cost and lead to over-refinement of subdomains with less complex solution behavior.

Another challenge in BEM formulations for multilayered media is the strong coupling between interfaces, which causes errors on one boundary to propagate to others.  For computational efficiency, it is therefore crucial to achieve a comparable accuracy on all interfaces: over-resolving one boundary while others remain coarse does not improve the global error and only increases the computational cost.   The adaptive algorithm employed in this work addresses this issue by iteratively refining the boundary meshes based on interpolation error indicators. In this way, the system automatically identifies where refinement is needed, leading to a balanced distribution of error across interfaces and an overall more efficient discretization.
  
Our automatic mesh adaptation algorithm for multilayer problems has been shown to be optimal for simple configurations where analytical solutions are available. However, our experiments also reveal a critical point: the need for highly accurate evaluation of nearly singular integrals. These integrals receive relatively little attention outside the BEM community as the solution is usually not required with high precision near boundaries.  For multilayered media, particularly when we may want to model metamaterials where high contrasts are prevalent, it is important to address but they are often treated as post-processing options.  In the context of mesh adaptation, however, inaccuracies in these integrals can dominate the local error indicators near interfaces, misleading the refinement process and causing unnecessary mesh enrichment driven by quadrature errors rather than by the solution dynamics. Several methods for handling near-singular integrals exist in the literature, but we have not yet incorporated them into our implementation.

Though the results discussed here are in 2D, we expect that for multilayered media in 3D the reduction of computational costs would be more pronounced and the potential benefit even greater.  For BEM, our focus on the interfaces is natural, but this approach and our conclusions could be useful for other mesh-based methods for wave problems as well, including FEM or FD methods. To sum up, less is more:  the solution itself should dictate the discretization; it is not needed to refine everywhere for optimal accuracy.


\appendix
\section{Analytic Solution for Multilayer Transmission Problems}
\label{sec:analytic solutions}
The analytic solution for the multilayer transmission problem is available in the case of concentric, circular interfaces \cite{marrec_2004}.  In these cases, the boundary \(\Gamma_j\) for \(j = 0, \hdots, N - 1\), is parameterized by \(\mathbf{x}_j(\theta) = \left(R_j \cos\theta, R_j \sin\theta \right) \) with \(R_N < R_{N-1} < \hdots < R_1 < R_0\).  The incident field, for the choice $\mathbf{a} = (0, 1)$, may also be expressed in terms of \(r\) and \(\theta\) such that \(u^{in}(y) =  u^{in}(r, \theta)= e^{ik_0 r\sin \theta}\).

By the Jacobi-Anger expansion \cite{NIST:DLMF},
\begin{equation}
    e^{iz \sin\theta} = \sum_{m=0}^{\infty} J_m(z) e^{im\theta}
        \label{eq:jacobi-anger expansion}
\end{equation}
the expression of the incident field may be given by
\begin{equation*}
    u^{in}(r, \theta) =  e^{ik_0 r\sin\theta} = \sum_{m=0}^{\infty} J_m(k_0 r) e^{im\theta}
\end{equation*}
where \(J_m\) is the \(m\)th order Bessel function. In $\Omega_0$, the scattered field must satisfy both the Helmholtz equation \eqref{eq:Helmholtz Eq} and the Sommerfeld radiation condition \eqref{eq:sommerfeld}.  It takes the form of the following ansatz

\begin{equation}
    u^{sc}_0(r, \theta) \equiv \sum_{m=0}^{\infty} A_m H^{(1)}_m(k_0r)e^{im\theta}
\end{equation}
where $H^{(1)}_m$ are \(m\)th order Hankel functions of the first kind, which represent outgoing wave behavior, and $A_m \in \mathbb{C}$  are coefficients to determine.  Similarly, the behavior of interior fields can be captured by expansions in terms of \(J_m\) and \(H^{(1)}_m\) in the regions $\Omega_j$ for $j = 1, \hdots, N$ \cite{colton_helmholtz_2019}.  Altogether, we can construct the analytic expression for the total field solution \(u^{\text{ana}}_j \) in each subdomain $\Omega_j$ as:

\begin{equation}
\resizebox{0.8\textwidth}{!}{$
\begin{aligned}
u^{\text{ana}}_0(r, \theta) &= \sum_{m=0}^{\infty} \left[ A_m H_m^{(1)}(k_0 r) + J_m(k_0 r_0) \right] e^{i m \theta} && R_0 \leq r,  \quad \forall \theta,\\
u^{\text{ana}}_j(r, \theta) &= \sum_{m=0}^{\infty} \left[ B_{j,m} H_m^{(1)}(k_j r) + C_{j,m} J_m(k_j r) \right] e^{i m \theta} && R_{j-1} \leq r \leq R_j, \quad \forall \theta, \\
u^{\text{ana}}_N(r, \theta) &= \sum_{m=0}^{\infty} E_m J_m(k_N r) e^{i m \theta} &&  r \leq R_N, \quad \forall \theta,
\end{aligned}
$}
\label{eq:analytic solutions}
\end{equation}
with \(B_{j,m}\), \(C_{j,m}\), \(E_m \in  \mathbb{C}\).  By enforcing the transmission conditions \eqref{eq:transmission problem}  at each interface $\Gamma_j$ between regions \(\Omega_j\) and \(\Omega_{j+1}\), we yield a solvable linear system for the unknown coefficients $A_m$, $B_{j, m}$, $C_{j, m}$, and $E_m$.  The unknown coefficients are then used in \eqref{eq:analytic solutions} to determine the solution $u^{ana}$ anywhere in the domain.

\section{Analytic Solutions for the Sound-Hard Boundary Condition}
\label{app:analytic solutions.}
We consider here the sound-hard scattering problem. For the one-boundary problem, the solutions \(u_0^{ana}\)  and \(u_1^{ana}\) must satisfy the Helmholtz equation as defined in \eqref{eq:Helmholtz Eq} in the regions \(\Omega_0 \cup\Omega_1\).  The sound-hard boundary condition on the interface \(\Gamma_0\) defined at radius \(R_0\) is \( \frac{\partial u_0}{\partial n_0} = 0\).

As similarly stated in \eqref{eq:analytic solutions}, the solutions in \(\Omega_0\) and \(\Omega_1\) respectively can be expressed as:
\[
\left\{
\begin{aligned}
    u_0^{ana}(r, \theta) &= \sum_{m=0}^\infty \left[ A_m H^{(1)}_m(k_0 r) + J_m(k_0 r)\right] e^{im\theta} \\
    u_1^{ana}(r, \theta) &= 0
\end{aligned}
\right.
\]
The coefficient \(A_m\) must satisfy:

\begin{equation*}
A_m k_0 H^{'(1)}_m(k_0R_0) + k_0 J'_m(k_0 R_0) = 0 .
\end{equation*}

Thus the coefficient \(A_m\) can be found as:

\begin{equation*}
    A_m =  \frac{-J'_m(k_0R_0)}{H^{'(1)}_m(k_0 R_0)}\\ 
\end{equation*}
Finally, the analytic solution for the sound-hard boundary case for a single boundary is given as:
\[
\left\{
\begin{aligned}
    u^{ana}_0(r, \theta) &= \sum_{m=0}^\infty \left[ J_m(k_0r) 
    - \frac{J'_m(k_0 R_0)}{H^{'(1)}_m(k_0R_0)} \, H^{(1)}_m(k_0 r)\right] e^{im\theta} \\
    u^{ana}_1(r, \theta) &= 0
\end{aligned}
\right.
\]

\bibliographystyle{elsarticle-num}

\bibliography{main}

\end{document}